\documentclass[12pt]{amsart}
\usepackage{amscd}
\usepackage[mathscr]{eucal}
\usepackage{amssymb}
\usepackage{latexsym}
\usepackage{amsthm}
\usepackage{amsmath}
\usepackage{graphicx}

\theoremstyle{plain}

\newtheorem{thm}{Theorem}[section]
\newtheorem{cor}{Corollary}[section]
\newtheorem{lem}{Lemma}[section]
\newtheorem{prop}{Proposition}[section]

\theoremstyle{definition}

\newtheorem*{dfn}{Definition}
\newtheorem{prf}{Proof}[section]

\DeclareMathOperator*{\essinf}{ess\, inf}
\DeclareMathOperator*{\esssup}{ess\, sup}

\title{Convolution invariant linear functionals and applications to summability methods}
\author{Ryoichi Kunisada}

\address{Faculty of Education and Integrated Arts and Science, Waseda University, Shinjuku-ku, Tokyo 169-8050, Japan}

\email{rkunisada@aoni.waseda.jp}

\thanks{This paper is part of the research performed under a Waseda University Grant for Special Research Projects (Project number: 2019C-135).}


\date{}

\begin{document}
\maketitle

\begin{abstract}
We study topologically invariant means on $L^{\infty}(\mathbb{R})$, the set of all essentially bounded functions on the real line, and prove that invariance with respect to a single convolution operator is sufficient for a mean to be topologically invariant. We also consider some applications of this result to summability methods. In particular, the notion of almost convergence is introduced for a function in $L^{\infty}(\mathbb{R})$, and a Tauberian theorem concerning almost convergence and a summability method defined by a Wiener kernel is obtained. Further, for the $C_{\infty}$ summability method, which is defined by the limit of H\"{o}lder summability methods, we provide a necessary and sufficient condition for a given function to be $C_{\infty}$ summable. 
\end{abstract}

\bigskip

\section{Introduction}
The Amenability of locally compact groups is of great importance for many branches of mathematics, including representation theory, operator algebras, and group theory (\cite{Pat}). Herein, since we are mainly interested in the application of such amenability to summability methods, we exclusively work with the special groups $\mathbb{R}$ and $\mathbb{R}^{\times}=(0, \infty)$, which are the additive group of the real field and the positive multiplicative group of $\mathbb{R}$, respectively. However, we remark that most of our results are valid for general locally compact abelian groups.

We give several definitions for the amenability of $\mathbb{R}$. 
Let $L^{\infty}(\mathbb{R})$ be the set of all essentially bounded functions on $\mathbb{R}$, and
let $L^1(\mathbb{R})$ be the group algebra of $\mathbb{R}$. Further, let $\mathscr{L}^1_{\flat}(\mathbb{R})$ be the set of functions $f \in L^1(\mathbb{R})$ such that $f \ge 0$ and $\hat{f}(0) = 1$. In the present study, we exclusively consider real-valued functions.

Let us define the action of $\mathbb{R}$ on $L^{\infty}(\mathbb{R})$ by  
\[T_s : L^{\infty}(\mathbb{R}) \rightarrow L^{\infty}(\mathbb{R}), \quad (T_s\phi)(x) = \phi(x+s) \ (s \in \mathbb{R}), \]
where $s \in \mathbb{R}$ and $\phi \in L^{\infty}(\mathbb{R})$.
For simplicity, for each $s \in \mathbb{R}$, we occasionally denote by $f_s(x)$ the translate $f(x+s)$ of a function $f(x)$ of a variable $x \in \mathbb{R}$. We further define its dual action on $L^{\infty}(\mathbb{R})^*$, the dual space of $L^{\infty}(\mathbb{R})$,  by $T_s^*\varphi(\phi) = \varphi(T_s\phi)$, where $\varphi \in L^{\infty}(\mathbb{R})^*$. For $\varphi$ in $L^{\infty}(\mathbb{R})^*$, $\varphi$ is said to be a {\it mean} on $L^{\infty}(\mathbb{R})$ if $\varphi$ is positive, i.e., $\varphi(\phi) \ge 0$ if $\phi \ge 0$, and $\|\varphi\|=1$. Note that this is equivalent to the condition that $\varphi(1) = \|\varphi\| =1$. A mean $\varphi$ on $L^{\infty}(\mathbb{R})$ is called an {\it invariant mean} if $T^*_s\varphi=\varphi$ holds for every $s \in \mathbb{R}$. Let $\mathscr{I}(\mathbb{R})$ be the set of all invariant means on $L^{\infty}(\mathbb{R})$. It is well known that the additive group $\mathbb{R}$ is amenable, that is, there exist invariant means on $L^{\infty}(\mathbb{R})$.

Next, we introduce another notion concerning amenability. Let us define an action of $L^1(\mathbb{R})$ on $L^{\infty}(\mathbb{R})$. For any $f \in L^1(\mathbb{R})$, the symbol $F$ denotes the convolution operator on $L^{\infty}(\mathbb{R})$ defined as follows:
\[F : L^{\infty}(\mathbb{R}) \rightarrow L^{\infty}(\mathbb{R}), \quad (F\phi)(x) = (f * \phi)(x), \quad \phi \in L^{\infty}(\mathbb{R}), \]
where $*$ is the convolution defined by
\[(f * \phi)(x) = \int_{-\infty}^{\infty} \phi(x-t)f(t)dt = \int_{-\infty}^{\infty} \phi(t)f(x-t)dt, \quad x \in \mathbb{R}. \]
Then, with the mapping $L^1(\mathbb{R}) \times L^{\infty}(\mathbb{R}) \rightarrow L^{\infty}(\mathbb{R})$, $(f, \phi) \mapsto F\phi$, $L^{\infty}(\mathbb{R})$ becomes an $L^1(\mathbb{R})$-modulo. We also consider the dual action on $L^{\infty}(\mathbb{R})^*$ defined by $F^*\varphi(\phi) = \varphi(F\phi)$, where $\varphi \in L^{\infty}(\mathbb{R})$. The following notion is closely related to our main objective.

\begin{dfn}
A mean $\varphi$ on $L^{\infty}(\mathbb{R})$ is said to be {\it topologically invariant} if 
\[F^*\varphi = \varphi \]
for every $f \in \mathscr{L}^1_{\flat}(\mathbb{R})$.
\end{dfn}
Let $\mathscr{T}(\mathbb{R})$ be the set of all topologically invariant means on $L^{\infty}(\mathbb{R})$. We know that for $\varphi \in \mathscr{T}(\mathbb{R})$, $\varphi(\phi_s) = \varphi(f*\phi_s) = \varphi(f_s*\phi) = \varphi(\phi)$. Hence, $\mathscr{T}(\mathbb{R}) \subseteq \mathscr{I}(\mathbb{R})$ immediately follows.
It is also straightforward to establish that for discrete amenable groups, such as the additive group of integers $\mathbb{Z}$, the invariant means and topologically invariant means coincide. However, in contrast to discrete groups, it is known that $\mathscr{T}(\mathbb{R}) \subsetneq \mathscr{I}(\mathbb{R})$ (see \cite{Pat}). Topologically invariant means are more manageable than simple invariant means and well suited to harmonic analysis.

Our main objective of this paper is means on $L^{\infty}(\mathbb{R})$ invariant with respect to a single convolution operator, which turn out to be equal to topologically invariant means. That is, means $\varphi$ on $L^{\infty}(\mathbb{R})$ for which $F^*\varphi = \varphi$ for a fixed $F$. Let us denote the set of all $F$-invariant linear functionals on $L^{\infty}(\mathbb{R})$ by $M_F$. $M_F$ is obviously a closed subspace of $L^{\infty}(\mathbb{R})^*$. We denote by $\mathscr{M}_F$ the subset of $M_F$ that consists of all $F$-invariant means. Note that $\mathscr{M}_F$ is a weak*-compact convex subset of $L^{\infty}(\mathbb{R})^*$. By definition, if $f \in \mathscr{L}^1_{\flat}(\mathbb{R})$ then $\mathscr{T}(\mathbb{R}) \subseteq \mathscr{M}_F$ holds true since $\mathscr{T}(\mathbb{R}) = \cap_{f \in \mathscr{L}^1_{\flat}(\mathbb{R})} \mathscr{M}_F$. 

Note that our study includes, as a special case, the invariant means with respect to the Hardy operator, which has been studied in the literature (\cite{Car, Kuni, Suk}). We provide simpler proofs of some of the results in these papers.

The property of a convolution operator $F$ depends on the behavior of the Fourier transform $\hat{f}$ of $f$, where $\hat{f}$ is defined by
\[\hat{f}(\xi) = \int_{-\infty}^{\infty} f(x) e^{-i\xi x} dx, \quad \xi \in \mathbb{R}. \]
Let $\mathscr{L}^1_*(\mathbb{R})$ be the set of functions $f$ in $L^1(\mathbb{R})$ such that $\hat{f}(\xi) = 1$ only at the point $\xi = 0$.
Further, we define $\mathscr{L}^1_{\sharp}(\mathbb{R})$ as the set of functions $f$ in $L^1(\mathbb{R})$ such that $f \ge 0$, $\hat{f}(0) = 1$, and $|\hat{f}(\xi)| = 1$ only at the point $\xi = 0$. 
By definition, it is clear that $\mathscr{L}^1_{\sharp}(\mathbb{R}) \subseteq \mathscr{L}^1_*(\mathbb{R})$ and $\mathscr{L}^1_{\sharp}(\mathbb{R}) \subseteq \mathscr{L}^1_{\flat}(\mathbb{R})$. We will show in the next section that, in fact, $\mathscr{L}^1_{\flat}(\mathbb{R}) = \mathscr{L}^1_{\sharp}(\mathbb{R})$ holds. Thus, it follows that $f \in \mathscr{L}^1_{\flat}(\mathbb{R})$ if and only if $f \in \mathscr{L}^1_*(\mathbb{R})$ and $f \ge 0$.

The first main result of this study is an analytic expression for the sublinear functional 
\[\overline{p}_F(\phi) = \sup_{\varphi \in \mathscr{M}_F} \varphi(\phi), \quad \phi \in L^{\infty}(\mathbb{R}). \]
Note that the functional $\overline{p}_F$ yields the maximal value of the elements of $\mathscr{M}_F$ for each fixed $\phi \in L^{\infty}(\mathbb{R})$. As will be explained in Section 2, the functional $\overline{p}_F$ plays an important role in the study of $M_F$.
More precisely, our result reads as follows. Let us define the sublinear functional $\overline{P}$ on $L^{\infty}(\mathbb{R})$ by
\[\overline{P}(\phi) = \lim_{\theta \to \infty} \limsup_{x \to \infty} \frac{1}{\theta} \int_x^{x+\theta} \phi(t)dt, \quad \phi \in L^{\infty}(\mathbb{R}). \]
Then, for any $F$ induced by $f \in \mathscr{L}^1_*(\mathbb{R})$, we will show that
\[\overline{p}_F(\phi) = \overline{P}(\phi) \]
holds for every $\phi \in L^{\infty}(\mathbb{R})$. In particular, this result implies that each set of $F$-invariant means are equal to one another for convolution operators $F$ induced by $f \in \mathscr{L}^1_{\flat}(\mathbb{R})$; thus, $\mathscr{M}_F = \mathscr{T}(\mathbb{R})$ holds true. This means that the invariance with respect to a single convolution operator $F$ with $f \in \mathscr{L}^1_{\flat}(\mathbb{R})$ is sufficient for a mean $m$ on $L^{\infty}(\mathbb{R})$ to be topologically invariant. This result forms a basis of the remainder of this paper. Moreover, this result is considered to be a refined characterization of topologically invariant means on $L^{\infty}(\mathbb{R})$ and may be valuable for the theory of amenability.

Furthermore, for an operator $F$ induced by $f \in \mathscr{L}^1_{\flat}(\mathbb{R})$, we show that each $\varphi \in M_F$ can be expressed uniquely as
\[\varphi = \alpha\varphi_+ - \beta\varphi_-, \quad \|\varphi\| = \alpha\|\varphi_+\| + \beta\|\varphi_-\|, \]
where $\varphi_+, \varphi_- \in \mathscr{M}_F$ and $\alpha, \beta \ge 0$. In other words, the elements of $M_F$ admit the Jordan decomposition. 

We also consider an interesting representation of the sublinear functional $\overline{P}$ as an infinite iteration of a sublinear functional related to $F$. Let $\overline{F} : L^{\infty}(\mathbb{R}) \rightarrow \mathbb{R}$ be defined by
\[\overline{F}(\phi) = \limsup_{x \to \infty} (F\phi)(x) = \limsup_{x \to \infty} \int_{-\infty}^{\infty} \phi(x-t)f(t)dt, \quad \phi \in L^{\infty}(\mathbb{R}). \]
We now consider the sublinear functionals $\overline{F}_k, k =1,2, \ldots$ defined by the iteration of $F$. For $k \ge 1$, we define $\overline{F}_k : L^{\infty}(\mathbb{R}) \rightarrow \mathbb{R}$ inductively as follows:
\[\overline{F}_k(\phi) := \overline{F}_{k-1}(F\phi) = \limsup_{x \to \infty} (F^k\phi)(x).  \]
Note that one can write $\overline{F}_k(\phi) = \limsup_{x \to \infty} (f^{*k} * \phi)(x)$, where $f^{*k}$ is the k-th power of $f$ with respect to the convolution product.
For an operator $F$ induced by $f \in \mathscr{L}^1_{\flat}(\mathbb{R})$, we can easily verify that
\[\overline{F}(\phi) = \overline{F}_1(\phi) \ge \overline{F}_2(\phi) \ge \cdots \ge \overline{F}_k(\phi) \ge \cdots \]
for every $\phi \in L^{\infty}(\mathbb{R})$. Thus, we can define the sublinear functional $\overline{F}_{\infty} : L^{\infty}(\mathbb{R}) \rightarrow \mathbb{R}$ as
\[\overline{F}_{\infty}(\phi) = \lim_{k \to \infty} \overline{F}_k(\phi), \quad \phi \in L^{\infty}(\mathbb{R}). \]
Let $F$ be an operator induced by $f \in \mathscr{L}^1_{\flat}(\mathbb{R})$. Then, for every $\phi \in L^{\infty}(\mathbb{R})$, we have
\[\overline{F}_{\infty}(\phi) = \overline{P}(\phi), \]
which is the second main result of this study. From the viewpoint of summability methods, this result can be interpreted as a representation of the relationship between a summability method defined by a convolution operator and the summability method defined by topologically invariant means. As will be described in Section 5, the sublinear functionals $\overline{F}$ and $\overline{P}$ give rise to summability methods for functions on $\mathbb{R}$, which we call the $F$ summability method and $P$ summability method, respectively. In particular, the class of $F$ summability methods includes most of the classical methods, including the Ces\`{a}ro, Abel, and Lambert methods (see \cite{Kor}). We show that for a function $\phi$ in $L^{\infty}(\mathbb{R})$, the $F$ summability of $\phi$ implies the $P$ summability of $\phi$. We also provide a Tauberian condition under which the converse implication holds. This is the third main result of this study.

The remainder of this paper is organized as follows. In Section 2, we present some preliminary results from the theory of topological linear spaces and the Fourier analysis of the Banach algebra $L^1(\mathbb{R})$. A continuous analogue of Banach limits is introduced. Further, we establish elementary results concerning the Fourier transforms of elements in $\mathscr{L}^1_{\flat}(\mathbb{R})$. 

In Section 3, we prove the first main result. That is, we obtain a characterization of convolution invariant functionals. In Section 4, we prove the second main result concerning the infinite iteration of sublinear functionals induced by convolution operators. For this purpose, we use a theorem of Katznelson and Tzafriri from operator theory.

In Section 5, we present an application to summability methods. By applying the results in Section 4, we present a Tauberian theorem involving P and F summability methods. 

In Section 6, we deal with a multiplicative version of the results in Sections 3, 4, and 5. Herein, the proofs goes as the additive case. Accordingly, we occasionally present only results without their proofs.

In Section 7, we consider Ces\`{a}ro invariant functionals. They can be viewed as a discrete analogue of Hardy invariant functionals, from which similar results can be obtained thorough elementary arguments. We also present some results on the $C_{\infty}$ summability method.

\section{Preliminaries}
Since we are concerned with weak*-compact convex subsets of dual spaces of $L^{\infty}$-spaces, the following version of the Krein-Milman theorem plays an important role (see \cite{Bou}).

\begin{prop}
Let $X$ be a Banach space and $X^*$ be its dual space. Further, let $\mathcal{C}$ be a weak*-compact convex subset of $X^*$ and $S \subseteq \mathcal{C}$. The following assertions are then equivalent:
\begin{flushleft}
$(1)$ $\sup_{\varphi \in S} \varphi(x) = \sup_{\varphi \in \mathcal{C}} \varphi(x)$ holds for each $x \in X$.
\end{flushleft} 

\noindent
$(2)$ $\mathcal{C} = \overline{co}(S)$, where $\mathcal{C}$ is the closed convex hull of $S$.

\noindent
$(3)$ The closure $\overline{S}$ of $S$ contains all extreme points of $\mathcal{C}$.
\end{prop}

Let $(\Omega, \mathscr{F}, \mu)$ be a measure space and $L^{\infty}(\mu)$ be the set of essentially bounded functions on $\Omega$.
Further, let $\mathscr{C}$ be the set of all weak*-compact convex subsets of the positive part of the unit sphere $S^+_{L^{\infty}(\mu)^*} = \{\varphi \in L^{\infty}(\mu)^* : \varphi \ge 0, \|\varphi\|=1\}$ of $L^{\infty}(\mu)^*$, the dual space of $L^{\infty}(\mu)$.
Let $\mathscr{S}$ be the set of all sublinear functionals $q$ on $L^{\infty}(\mu)$ such that $q \ge 0$ and $q(1) = 1$. 
Consider the partially ordered sets $(\mathscr{C}, \subseteq)$ and $(\mathscr{S}, \le)$, where $\subseteq$ denotes the inclusion of subsets and $\le$ denotes the pointwise order of functionals. Based on the above proposition, we have the isomorphism between these partially ordered sets defined by
\[\mathscr{C} \ni \mathcal{C} \rightarrow q(x) = \sup_{\varphi \in \mathcal{C}} \varphi(x) \in \mathscr{S}.  \]
The inverse mapping is given by
\[\mathscr{S} \ni q \rightarrow \mathcal{C} := \{\varphi : \varphi(x) \le q(x) \ \text{for each} \ x \in L^{\infty}(\mu)\} \in \mathscr{C}. \]
This fact illustrates the importance of the study of the sublinear functional $\sup_{\varphi \in \mathcal{C}} \varphi(x)$ for a given class $\mathcal{C}$ of linear functionals since it represents the size of $\mathcal{C}$ in $S^+_{L^{\infty}(\mu)^*}$. Furthermore, owing to the equivalence of $(1)$ and $(3)$, it can be used to obtain some information about the extreme points $ex(\mathcal{C})$ of $\mathcal{C}$. Additionally, according to the equivalence of $(1)$ and $(2)$, one may obtain a set of functionals $A$ in $\mathcal{C}$ with a simple form such that they generate all the elements of $\mathcal{C}$ by taking its closed convex hull, i.e., $\mathcal{C} = \overline{co}(A)$. 

Note that one can obtain the minimal values of $\mathcal{C}$ from the maximal value functional $q(x) = \sup_{\varphi \in \mathcal{C}} \varphi(x)$ of $\mathcal{C}$. In fact, since $\varphi(-x) \le q(-x)$, we have $\underline{q}(x):= -\overline{q}(-x) \le \varphi(x)$. Hence, we obtain the range of $\varphi \in \mathcal{C}$ for each fixed $x \in L^{\infty}(\mu)$:
\[\underline{q}(x) \le \varphi(x) \le \overline{q}(x). \]
These inequalities are strict in the sense that for any real number $\alpha \in [\underline{q}(x), \overline{q}(x)]$, there exists a $\varphi \in \mathcal{C}$ such that $\varphi(x) = \alpha$. This is a consequence of the Hahn-Banach theorem.

In particular, for the set of means $\mathscr{M}(\mathbb{R})$ on $L^{\infty}(\mathbb{R})$, we obtain the following result.
\begin{prop}
For an element $n \in L^{\infty}(\mathbb{R})^*$, $n$ is in $\mathscr{M}(\mathbb{R})$ if and only if
\[n(\phi) \le \esssup_{x \in \mathbb{R}} \phi(x) \]
holds for every $\phi \in L^{\infty}(\mathbb{R})$. In particular, for a mean $m$ and $\phi \in L^{\infty}(\mathbb{R})$, we have
\[\essinf_{x \in \mathbb{R}} \phi(x) \le m(\phi) \le \esssup_{x \in \mathbb{R}} m(\phi). \]
\end{prop}

Let $C_{bu}(\mathbb{R})$ be the set of all bounded, uniformly continuous functions on $\mathbb{R}$ and $C_{bu}(\mathbb{R})^*$ be its dual space. Clearly, $C_{bu}(\mathbb{R})$ is a closed subalgebra of $L^{\infty}(\mathbb{R})$. We now study translation-invariant linear functionals on $C_{bu}(\mathbb{R})$, which play an important role in the study of convolution-invariant linear functionals on $L^{\infty}(\mathbb{R})$.

 Let $M_T$ be the set of all continuous linear functionals $\varphi$ on $C_{bu}(\mathbb{R})$ that are invariant under translations on $\mathbb{R}$, i.e., $\varphi \in C_{bu}(\mathbb{R})^*$, for which $T_s^*\varphi = \varphi$ holds for every $s \in \mathbb{R}$, where $T_s^*$ is the adjoint operator of $T_s$. Let $\mathscr{M}_T$ be the set of invariant means, i.e., the subset of $M_T$, the elements $\varphi$ of which satisfy the conditions $\varphi \ge 0$ and $\|\varphi\|=1$. $\mathscr{M}_T$ is a weak*-compact convex subset of $C_{bu}(\mathbb{R})^*$. We then obtain the following result readily from Banach lattice theory.

\begin{thm}
Let $\varphi \in M_T$. There exist $\varphi_+$ and $\varphi_-$ in $\mathscr{M}_T$ such that 
\[\varphi = \alpha \varphi_+ - \beta \varphi_-, \quad \|\varphi\| = \alpha\|\varphi_+\| + \beta\|\varphi_-\| \]
holds for some constants $\alpha, \beta \ge 0$.
\end{thm}

Hence, $M_T$ is generated by $\mathscr{M}_T$. Thus, it is sufficient to consider $\mathscr{M}_T$ for the study of $M_T$. The following result provides a necessary and sufficient condition that $\varphi \in C_{bu}(\mathbb{R})^*$ belongs to $\mathscr{M}_T$ (see \cite{Kuni} for the proof). Let the sublinear functional $\overline{P} : L^{\infty}(\mathbb{R}) \rightarrow \mathbb{R}$ be that defined in the introduction.

\begin{thm}
For $\varphi \in C_{bu}(\mathbb{R})^*$, $\varphi \in \mathscr{M}_T$ if and only if
\[\varphi(\phi) \le \overline{P}(\phi) \]
holds for every $\phi \in C_{bu}(\mathbb{R})$.
\end{thm}

We remark that the functional $\overline{P}$ is equal to the following one:
\[\overline{P}_1(\phi) = \lim_{\theta \to \infty} \sup_{x \ge 0} \frac{1}{\theta} \int_x^{x+\theta} \phi(t)dt, \quad \phi \in L^{\infty}(\mathbb{R}), \]
where $\limsup_{x \to \infty}$ is replaced by $\sup_{x \ge 0}$ in the definition $\overline{P}$. This functional was adopted in \cite{Suk}, where an assertion equivalent to Theorem 2.2 was proved.

We mention that the class $\mathscr{M}_T$ of linear functionals can be viewed as a continuous analogue of the classical notion of Banach limits. Recall that a Banach limit is a continuous linear functional $\varphi$ on $l_{\infty}$ such that $\varphi \ge 0$, $\|\varphi\| = 1$, and $\varphi$ invariant with respect to the translation operator $T: l_{\infty} \rightarrow l_{\infty}$, $T\phi(n) = \phi(n+1)$. Let $\overline{B}$ be the sublinear functional on $l_{\infty}$ defined by
\[\overline{B}(\phi) = \lim_{k \to \infty} \limsup_{n \to \infty} \frac{1}{k} \sum_{i=n}^{n+k-1} \phi(i), \quad \phi \in l_{\infty}. \]
The following result was reported in \cite{Jer}: for $\varphi \in l_{\infty}^*$, $\varphi$ is a Banach limit if and only if $\varphi(\phi) \le \overline{B}(\phi)$ holds for every $\phi \in l_{\infty}$ (see also \cite{Suc}). Notice that the functional $\overline{P}$ is a continuous analogue of $\overline{B}$ obtained by replacing the discrete summation with integration. For a more detailed exposition of translation-invariant linear functionals on $C_{bu}(\mathbb{R})$ (or $L^{\infty}(\mathbb{R})$), see \cite{Car, Kuni, Suk}.

We use some results from the theory of Fourier analysis on $L^1(\mathbb{R})$ and $L^{\infty}(\mathbb{R})$. We refer the reader to \cite{Rud} for details. Recall that the group algebra $L^1(\mathbb{R})$ is a Banach algebra, the product of which is the convolution $*$.
Let $I$ be a closed ideal of $L^1(\mathbb{R})$. The zero set $Z(I)$ of an ideal $I$ is defined by $Z(I) = \{\xi \in \mathbb{R} : \hat{f}(\xi) = 0 \ \forall f \in I\}$. $Z(I)$ is always closed and for each closed set $E$ of $\mathbb{R}$, there exists a closed ideal $I$ such that $Z(I) = E$. In particular, a closed set $E$, which is the zero set of a unique ideal $I$ of $L^1(\mathbb{R})$, is referred to as a {\it spectral synthesis set}. In other words, $E \subset \mathbb{R}$ is a spectral synthesis set if and only if $I = J$ holds for any closed ideals $I$ and $J$ of $L^1(\mathbb{R})$ with $Z(I) = Z(J) = E$. The following is a sufficient condition for a closed set $E$ of $\mathbb{R}$ to be a spectral synthesis set (see \cite{Rud}).
\begin{thm}
Let $E$ be a closed set of $\mathbb{R}$, the boundary of which contains no perfect set. Then, $E$ is a spectral synthesis set.
\end{thm}
Note that the celebrated Wiener's Tauberian theorem can be viewed as a special case of this result for $E = \emptyset$.
\begin{thm}
Let $I$ be a closed ideal of $L^1(\mathbb{R})$ and $Z(I) = \emptyset$. Then, $I = L^1(\mathbb{R})$ holds.
\end{thm}
In this paper, we also need a special case of this result in which $E$ is the singleton $\{0\}$.

Now, recall that $L^{\infty}(\mathbb{R})$ is a dual space of $L^1(\mathbb{R})$. For any closed ideal $I$ of $L^1(\mathbb{R})$, its annihilator $\Phi = I^{\bot} = \{\phi \in L^{\infty}(\mathbb{R}) : \langle f, \phi \rangle =0$ for every $f \in I \}$ in $L^{\infty}(\mathbb{R})$ is a weak*-closed translation invariant subspace. Conversely, the annihilator $I = \Phi^{\bot}$ of each weak*-closed translation subspace $\Phi$ is a closed ideal of $L^1(\mathbb{R})$. By using the duality between these  classes of subspaces of $L^1(\mathbb{R})$ and $L^{\infty}(\mathbb{R})$, the above result on closed ideals of $L^1(\mathbb{R})$ can be transferred to the context of weak*-closed invariant subspaces of $L^{\infty}(\mathbb{R})$, which play an important role in our study.
For a weak*-closed invariant subspace $\Phi$ of $L^{\infty}(\mathbb{R})$, its spectrum $\sigma(\Phi)$ in the sense of spectral synthesis is defined by $\sigma(\Phi) = \{\lambda \in \mathbb{C} : e^{i{\lambda}x} \in \Phi\}$, i.e., the continuous characters of $\mathbb{R}$ contained in $\Phi$. $\sigma(\Phi)$ is always closed, and for each closed set $E$ of $\mathbb{R}$, there exists a weak*-closed translation invariant subspace of $L^{\infty}(\mathbb{R})$ such that  $\sigma(\Phi) = E$. A closed subset $E$, which is the spectrum $\sigma(\Phi)$ of the unique $\Phi$, is called a spectral synthesis set. This definition of spectral synthesis sets is equivalent to the aforementioned one, as can be verified by the relation $\sigma(\Phi) = Z(I)$ for $I = \Phi^{\bot}$. Suppose that $\sigma(\Phi)$ is a spectral synthesis set. Let us consider the weak*-closed translation-invariant subspace $\Phi_1$ generated by $\sigma(\Phi)$. Then, $\Phi = \Phi_1$ holds. These facts, in conjunction with Theorem 2.3, imply the following result.
\begin{thm}
Let $\Phi$ be a weak*-closed translation-invariant subspace of $L^{\infty}(\mathbb{R})$ such that the boundary of its spectrum $\sigma({\Phi})$ contains no perfect set. Then, $\Phi$ is equal to the weak*-closed translation-invariant subspace generated by $\sigma(\Phi)$.
\end{thm}

In particular, we need the following special case of the above theorem, where the spectrum is the singleton $\{0\}$.
\begin{cor}
Let $\Phi$ be a weak*-closed invariant subspace of $L^{\infty}(\mathbb{R})$ with $\sigma(\Phi) = \{0\}$. Then, $\Phi$ is the subspace of $L^{\infty}(\mathbb{R})$ consisting of the constant functions.
\end{cor}

Now, we prove the assertion that $\mathscr{L}^1_{\flat}(\mathbb{R}) = \mathscr{L}^1_{\sharp}(\mathbb{R})$.
\begin{thm}
$\mathscr{L}^1_{\flat}(\mathbb{R}) = \mathscr{L}^1_{\sharp}(\mathbb{R})$ holds.
\end{thm}

\begin{prf}
It is known by definition that $\mathscr{L}^1_{\sharp}(\mathbb{R}) \subseteq \mathscr{L}^1_{\flat}(\mathbb{R})$. Thus, it is sufficient to show the opposite inclusion. Suppose that $f \in \mathscr{L}^1_{\flat}(\mathbb{R})$.
First, we show that $\hat{f}(\xi) = 1$ only at the point $\xi = 0$. In fact, suppose that for a $\xi \in \mathbb{R} \setminus \{0\}$,
\[\hat{f}(\xi) = \int_{-\infty}^{\infty} f(x)e^{-i \xi x} dx = \int_{-\infty}^{\infty} f(x) \cos(\xi x)dx - i\int_{-\infty}^{\infty} f(x)\sin(\xi x)dx = 1. \]
It follows immediately that the first term is $1$ and the second term is $0$. From the assumption that $\hat{f}(0) = 1$ and $f \ge 0$, if 
\[\int_{-\infty}^{\infty} f(x)\cos(\xi x)dx = 1 \]
holds, then it is necessary that $f = 0$ on the subset of $\mathbb{R}$ at which $\cos (\xi x)$ is negative, i.e., $\{x \in \mathbb{R} : \cos (\xi x) < 0\} = \cup_{n = -\infty}^{\infty} (\frac{\pi}{2\xi}+\frac{2n\pi}{\xi}, \frac{3}{2\xi}\pi+\frac{2n\pi}{\xi})$. Assume that this condition is satisfied. Then we have
\begin{align}
\int_{-\infty}^{\infty} f(x)\cos (\xi x)dx &= \sum_{n=-\infty}^{\infty} \int_{-\frac{\pi}{2\xi}+\frac{2n\pi}{\xi}}^{\frac{\pi}{2\xi}+\frac{2n\pi}{\xi}} f(x) \cos(\xi x)dx  \notag \\
&< \sum_{n=-\infty}^{\infty}  \int_{-\frac{\pi}{2\xi}+\frac{2n\pi}{\xi}}^{\frac{\pi}{2\xi}+\frac{2n\pi}{\xi}} f(x)dx = 1\notag
\end{align}
since the Lebesgue measure of the set of points at which $\cos (\xi x) = 1$ is $0$. This contradicts the assumption.
Since the space $\mathscr{L}^1_{\flat}(\mathbb{R})$ is translation invariant, we can generalize the above result to the case of absolute values of $\hat{f}(\xi)$ through the following argument.
Suppose that for a $\xi_0 \in \mathbb{R} \setminus \{0\}$,
\[|\hat{f}(\xi_0)| = \left|\int_{-\infty}^{\infty} f(x)e^{-i\xi_0 x}dx \right| = 1. \]
Let the argument of $\hat{f}(\xi_0)$ be $\theta_0$. Hence, we have 
\begin{align}
e^{-i \theta_0} \hat{f}(\xi_0) = 1 &\Longleftrightarrow e^{-i \theta_0} \int_{-\infty}^{\infty} f(x)e^{-i \xi_0 x}dx = 1 \notag \\ 
&\Longleftrightarrow \int_{-\infty}^{\infty} f(x)e^{-i(\xi_0 x + \theta_0)}dx = 1. \notag 
\end{align}
By using $x = y - \theta_0/\xi_0$ and integration by substitution, we have
\[\int_{-\infty}^{\infty} f(x)e^{-i(\xi_0 x + \theta_0)}dx = \int_{-\infty}^{\infty} f\left(y - \frac{\theta_0}{\xi_0}\right) e^{-i \xi_0 y}dy = 1. \]
This means that $(f_{-\theta_0/\xi_0}) \ \hat{} \ (\xi_0) = 1$. However, it is obvious that $f_{-\theta_0/\xi_0} \in \mathscr{L}^1_{\flat}(\mathbb{R})$. Since we have assumed $\xi_0 \not= 0$, this result contradicts the result presented above. This completes the proof.
\end{prf}

\section{Convolution-invariant functionals on $L^{\infty}(\mathbb{R})$}
In this section, we characterize $F$-invariant functionals, where the convolution operator $F$ is induced by the elements $f$ in $\mathscr{L}^1_*(\mathbb{R})$. First, we obtain a characterization of $F$-invariant functionals on $C_{bu}(\mathbb{R})$, which in turn is used to derive that of $F$-invariant functionals on $L^{\infty}(\mathbb{R})$. 
\begin{thm}
For any $\phi \in C_{bu}(\mathbb{R})$ and $\varphi \in L^{\infty}(\mathbb{R})^*$, 
\[(F^*\varphi)(\phi) = \int_{-\infty}^{\infty} \varphi(\phi_t)f(-t)dt \]
holds true.
\end{thm}

\begin{prf}
Since the function $T_s\phi : (-\infty, \infty) \rightarrow C_{bu}(\mathbb{R})$ is continuous, it is $C_{bu}(\mathbb{R})$-valued Bochner $f(x)dx$-integrable. Thus, from Corollary 2 of [20, p.134], we have
\begin{align}
(F^*\varphi)(\phi) &= \varphi(F\phi) = \varphi\left(\int_{-\infty}^{\infty} \phi(x-t)f(t)dt\right)  \notag \\
&= \varphi\left(\int_{-\infty}^{\infty} \phi_{-t}(x)f(t)dt\right) = \int_{-\infty}^{\infty} \varphi(\phi_{-t})f(t)dt. \notag
\end{align}
The proof is thus complete.
\end{prf}

The following corollary follows immediately.
\begin{cor}
For any $\phi \in C_{bu}(\mathbb{R})$ and $\varphi \in L^{\infty}(\mathbb{R})^*$, 
let $\psi(s) = \varphi(\phi_s)$, where $s \in \mathbb{R}$. Then, we have
\[(F^*\varphi)(\phi_s) = \int_{-\infty}^{\infty} \varphi(\phi_t)f(s-t)dt  = (\psi * f)(s). \]
\end{cor}

Note that $\psi$ is in $C_{bu}(\mathbb{R})$.
Hence, if $\varphi \in L^{\infty}(\mathbb{R})^*$ is $F$-invariant, then for any $\phi \in C_{bu}(\mathbb{R})$, we have
\[\psi(x) = F\psi(x), \]
where $\psi(x) = \varphi(\phi_x)$ with $x \in \mathbb{R}$. This implies that $\psi$ is an eigenfunction of the convolution operator $F$ with an eigenvalue one.  

Now, we study some basic properties of convolution operators $F$ on $L^{\infty}(\mathbb{R})$: 
\[F : L^{\infty}(\mathbb{R}) \rightarrow L^{\infty}(\mathbb{R}), \quad (F\phi)(x) = (f * \phi)(x). \] 
First, we show the weak* continuity of convolution operators.
\begin{thm}
$F$ is a weak*-continuous linear operator on $L^{\infty}(\mathbb{R})$.
\end{thm}

\begin{prf}
Suppose that $\phi_{\alpha}, \phi \in L^{\infty}(\mathbb{R})$, and $w^*\mathchar`-\lim_{\alpha} \phi_{\alpha} = \phi$, i.e., the net $\{\phi_{\alpha}\}$ converges to $\phi$ in the weak* sense. It is sufficient to show that $w^*\mathchar`-\lim_{\alpha} f * \phi_{\alpha} = f * \phi$. In other words, for every $g \in L^1(\mathbb{R})$, we show that
\[\lim_{\alpha} \int_{-\infty}^{\infty} (\phi_{\alpha} *f)(x)g(-x)dx = \int_{-\infty}^{\infty} (\phi * f)(x)g(-x)dx. \]
Notice that this equality is equal to
\[\lim_{\alpha} \int_{-\infty}^{\infty} \phi_{\alpha}(x)(f * g)(-x)dx = \int_{-\infty}^{\infty} \phi(x)(f * g)(-x)dx. \]
From the assumption that $\phi_{\alpha}$ converges $\phi$ in the weak* sense, we obtain the theorem.
\end{prf}
We can now determine the spectrum of $F$. From \cite{Jor} of Theorem 13.2, we have the following result.
\begin{thm}
For any convolution operator $F$, its spectrum $\sigma(F)$ is $\{0\} \cup \{\lambda \in \mathbb{C} : \hat{f}(x) = \lambda$ for some $x \in \mathbb{R}\}$. In particular, $\sigma(F) \setminus \{0\}$ consists of point spectra of $F$.
\end{thm}

We now consider the space of eigenfunctions of a convolution operator $F$. Let $\lambda \in \sigma_p(F)$ be an eigenvalue of $F$ and $E_{\lambda}(F) = \{\phi \in L^{\infty}(\mathbb{R}) : F\phi = \lambda\phi\}$ be the eigenfunction space of $F$ with respect to the eigenvalue $\lambda$. The following results concerning the space $E_{\lambda}(F)$ are now obtained.
\begin{thm}
$E_{\lambda}(F)$ is a weak*-closed invariant subspace of $L^{\infty}(\mathbb{R})$.
\end{thm}

\begin{prf}
Suppose that $\phi \in E_{\lambda}(F)$, i.e., $(f * \phi)(x) = \lambda\phi(x)$ for every $x \in \mathbb{R}$. Then, we have $(f * \phi_s)(x) = (f * \phi)(x+s) = \lambda\phi(x+s) =\lambda\phi_s(x)$ for every $s \in \mathbb{R}$. Thus, $E_{\lambda}(F)$ is translation invariant. The assertion that $E_{\lambda}(F)$ is weak*-closed is clear from Theorem 3.2.
\end{prf}

\begin{thm}
The spectrum of $E_{\lambda}(F)$ in the sense of spectral synthesis is $\{\xi \in \mathbb{R} : \hat{f}(\xi) = \lambda\}$.
\end{thm}

\begin{prf}
For any $\xi \in \mathbb{R}$, $e^{i\xi x} \in E_{\lambda}(F)$ if and only if $f * e_{\xi} = \lambda e_{\xi}$. This implies that $\hat{f}(\xi)e_{\xi} = \lambda e_{\xi}$, from which we obtain the result immediately.
\end{prf}

Based on these results, we can deduce the following characterization of $F$-invariant functionals on $C_{bu}(\mathbb{R})$. 
\begin{thm}
Let $f \in \mathscr{L}^1_*(\mathbb{R})$ and $\varphi \in C_{bu}(\mathbb{R})^*$. Then, $\varphi$ is $F$-invariant if and only if $\varphi$ is translation invariant. 
\end{thm}

\begin{prf}
Let us assume that $\varphi$ is $F$-invariant. We show that $\varphi(T_s\phi) = \varphi(\phi)$ for every $\phi \in C_{bu}(\mathbb{R})$ and $s \in \mathbb{R}$. From Corollary 3.1, for any $\phi \in C_{bu}(\mathbb{R})$, we have
\[\varphi(\phi_s) = \int_{-\infty}^{\infty} \varphi(\phi_t)f(s-t)dt. \]
Thus, we obtain $\psi(s) = \varphi(\phi_s) \in E_1(F)$. However, $\sigma(E_1(F)) = \{0\}$ from Theorem 3.5 and the assumption on $f$. This implies that $E_1(F) = \mathbb{R}$ from Corollary 2.1, i.e., that they are constant functions. Hence, we have $\varphi(\phi_s) = \psi(s) = \psi(0) = \varphi(\phi)$ for every $\phi \in C_{bu}(\mathbb{R})$ and $s \in \mathbb{R}$. Hence, $\varphi$ is translation invariant on $C_{bu}(\mathbb{R})$. 

Suppose that $\varphi$ is in $M_T$, that is, translation invariant. We show that $\varphi$ is $F$-invariant. As indicated in Theorem 3.1, we have
\begin{align}
(F^*\varphi)(\phi) &= \int_{-\infty}^{\infty} \varphi(\phi_t)f(-t)dt = \int_{-\infty}^{\infty} \varphi(\phi)f(-t) \notag \\
&= \varphi(\phi) \int_{-\infty}^{\infty} f(-t)dt = \varphi(\phi), \notag
\end{align}
which shows that $\varphi$ is F-invariant.
This completes the proof.
\end{prf}

The lemma below is necessary to obtain the characterization of $\mathcal{M}_F$.
\begin{lem}
For any $\phi \in L^{\infty}(\mathbb{R})$ and $f \in L^1(\mathbb{R})$ with $\hat{f}(0) = 1$, $\overline{P}(f * \phi - \phi) = 0$ holds.
\end{lem}

\begin{prf}
First, through direct computation, we can obtain
\begin{align}
\overline{P}(\phi - f * \phi) &= \lim_{\theta \to \infty} \limsup_{x \to \infty} \frac{1}{\theta} \int_x^{x + \theta} (\phi - f * \phi)(t)dt \notag \\
&= \lim_{\theta \to \infty} \limsup_{x \to \infty} \frac{1}{\theta} \int_x^{x + \theta} \left\{\phi(t) - \int_{-\infty}^{\infty} \phi(t-s)f(s)ds\right\}dt \notag \\
&= \lim_{\theta \to \infty} \limsup_{x \to \infty} \frac{1}{\theta} \int_x^{x + \theta} \int_{-\infty}^{\infty} \{\phi(t) - \phi(t-s)\}f(s)dsdt \notag \\
&= \lim_{\theta \to \infty} \limsup_{x \to \infty} \frac{1}{\theta} \int_x^{x + \theta} \{\phi(t) - \phi(t-s)\}dt \int_{-\infty}^{\infty} f(s)ds. \notag
\end{align}
For any $\varepsilon > 0$, we can choose $R > 0$ such that $\int_{-R}^R f(s)ds > 1-\varepsilon$. Observe that
\begin{align}
\left|\frac{1}{\theta} \int_x^{x + \theta} \{\phi(t) - \phi(t-s)\}dt \right| &= \left|\frac{1}{\theta} \int_{x-s}^x \phi(t)dt + \frac{1}{\theta} \int_{x+\theta-s}^{x+\theta} \phi(t)dt \right|  \notag \\
&\le \frac{2\|\phi\|_{\infty}s}{\theta} \notag
\end{align}
and
\[\left|\frac{1}{\theta} \int_x^{x + \theta} \{\phi(t) - \phi(t-s)\}dt \right| \le \frac{1}{\theta} \cdot 2\|\phi\|_{\infty}\theta = 2\|\phi\|_{\infty}. \]
Hence, we have 
\begin{align}
&\left|\frac{1}{\theta} \int_x^{x + \theta} \{\phi(t) - \phi(t-s)\}dt \int_{-\infty}^{\infty} f(s)ds\right| \notag \\
&\le \left|\frac{1}{\theta} \int_x^{x + \theta} \{\phi(t) - \phi(t-s)\}dt \int_{-R}^R f(s)ds \right| \notag \\
&+ \left| \frac{1}{\theta} \int_x^{x + \theta} \{\phi(t) - \phi(t-s)\}dt \int_{(-\infty, -R] \cup [R, \infty)} f(s)ds \right| \notag \\
&\le \frac{2\|\phi\|_{\infty}R}{\theta} + 2\|\phi\|_{\infty}\varepsilon, \notag
\end{align}
which tends to $0$ as $\theta$ tends to $\infty$. Thus, we obtain the lemma.
\end{prf}

\begin{cor}
For any $\phi \in L^{\infty}(\mathbb{R})$ and $f \in L^1(\mathbb{R})$ with $\hat{f}(0) =1$, $\overline{P}(f * \phi) = \overline{P}(\phi)$ holds.
\end{cor}

We provide a characterization of $F$-invariant means $\mathscr{M}_F$ for a convolution operator induced by an element in $\mathscr{L}^1_*(\mathbb{R})$.
\begin{thm}
Let $f \in \mathscr{L}^1_*(\mathbb{R})$ and $F$ be the induced convolution operator. For $\varphi \in L^{\infty}(\mathbb{R})^*$, $\varphi \in \mathscr{M}_F$ if  and only if 
\[\varphi(\phi) \le \overline{P}(\phi) \]
holds for every $\phi \in L^{\infty}(\mathbb{R})$.
\end{thm}

\begin{prf}
Suppose that $\varphi$ is in $\mathscr{M}_F$. From Theorem 3.6, $\varphi$ is translation invariant on $C_{bu}(\mathbb{R})$ and is thus in $\mathscr{M}_T$. Then, from Theorem 2.2,
\[\varphi(\phi) \le \overline{P}(\phi) \]
holds for every $\phi \in C_{bu}(\mathbb{R})$. For $\phi \in L^{\infty}(\mathbb{R})$, from Corollary 3.2 and the $F$-invariance of $\varphi$, we have 
\[\varphi(\phi) = \varphi(f * \phi) \le \overline{P}(f * \phi) = \overline{P}(\phi), \]
which proves the necessity. Conversely, if 
\[\varphi(\phi) \le \overline{P}(\phi) \]
holds for every $\phi \in L^{\infty}(\mathbb{R})$, then it is valid that
\[\underline{P}(\phi - f * \phi) \le \varphi(\phi - f * \phi) \le \overline{P}(\phi - f * \phi). \]
From the proof of Lemma 3.1, it is easy to see that $\underline{P}(\phi - f * \phi) = 0$ is also valid. Thus, we have
\[\varphi(\phi - f * \phi) = 0. \]
Hence, 
\[F^*\varphi(\phi) = \varphi(\phi) \]
holds for every $\phi \in L^{\infty}(\mathbb{R})$, which implies that $\varphi$ is $F$-invariant.
\end{prf}

\begin{cor}
Let $f \in \mathscr{L}_{\flat}(\mathbb{R})$. Then, $\mathscr{M}_F = \mathscr{T}(\mathbb{R})$ holds. In other words, a mean $m$ is $\mathscr{L}^1_{\flat}(\mathbb{R})$-invariant if and only if it is $F$-invariant.
\end{cor}

Moreover, under the additional assumption on $f$ that $f \ge 0$, we have the following theorem concerning general $F$-invariant functionals $M_F$.
\begin{thm}
Let $F$ be a convolution operator induced by an element $f$ of $\mathscr{L}^1_{\flat}(\mathbb{R})$. Then, for any $\varphi \in M_F$, $\varphi$ admits the Jordan decomposition in $M_F$, i.e., there exist some positive elements $\varphi_+$ and $\varphi_-$ in $M_F$ such that 
\[\varphi = \varphi_+ - \varphi_- , \quad \|\varphi\| = \|\varphi_+\| + \|\varphi_-\|    \]
hold.
\end{thm}

\begin{prf}
Let us denote by $\varphi_0$ the restriction of $\varphi$ to $C_{bu}(\mathbb{R})$. Since $\varphi_0$ is translation invariant on $C_{bu}(\mathbb{R})$ as per Theorem 3.6, $\varphi_0$ can be decomposed using Theorem 2.1 as
\[\varphi_0 = \varphi_{0, +} - \varphi_{0, -}, \]
where $\varphi_{0, +}, \varphi_{0, -} \in M_T$ are positive and $\|\varphi_0\| = \|\varphi_{0, +}\| +\|\varphi_{0, -}\|$. We define continuous linear functionals $\overline{\varphi}_{0, +}$ and $ \overline{\varphi}_{0, -}$ on $L^{\infty}(\mathbb{R})$, which are positive and $F$-invariant, as
\[\overline{\varphi}_{0, +}(\phi) = \varphi_{0, +}(F\phi), \quad \overline{\varphi}_{0, -}(\phi) = \varphi_{0, -}(F\phi). \]
The positivity of $\overline{\varphi}_{0, +}$ and $\overline{\varphi}_{0, -}$ follows from the positivity of $\varphi_{0,+}$, $\varphi_{0, -}$, and $F$. The $F$-invariance of $\overline{\varphi}_{0, +}$ and $\overline{\varphi}_{0, -}$ is straightforward from their definitions. Furthermore, note that both $\overline{\varphi}_{0, +}$ and $\overline{\varphi}_{0, -}$ are extensions of $\varphi_{0, +}$ and $\varphi_{0, -}$, respectively. In fact, since $\varphi_{0, +}$ and $\varphi_{0, -}$ are in $M_T$, they are also $F$-invariant on $C_{bu}(\mathbb{R})$. Thus, for any $\phi \in C_{bu}(\mathbb{R})$, we have $\overline{\varphi}_{0, +}(\phi) = \varphi_{0, +}(F\phi) = \varphi_{0, +}(\phi)$.
In the same way, we can show that $\overline{\varphi}_{0, -}(\phi) = \varphi_{0, -}(\phi)$ for every $\phi \in C_{bu}(\mathbb{R})$. Observe that
\[|\overline{\varphi}_{0, +}(\phi)| = |\varphi_{0, +}(f * \phi)| \le \|\varphi_{0, +}\| \cdot \|f * \phi\|_{\infty} \le \|\varphi_{0, +}\| \cdot \|\phi\|_{\infty}, \]
which implies that $\|\overline{\varphi}_{0, +}\| \le \|\varphi_{0, +}\|$. Note that $\|\overline{\varphi}_{0, +}\| \ge \|\varphi_{0, +}\|$ follows because $\overline{\varphi}_{0, +}$ is an extension of $\varphi_{0, +}$. Thus, we obtain $\|\overline{\varphi}_{0, +}\| = \|\varphi_{0, +}\|$. In the same way, we have $\|\overline{\varphi}_{0, -}\| = \|\varphi_{0, -}\|$.

We set $\overline{\varphi}_0 = \overline{\varphi}_{0, +} - \overline{\varphi}_{0, -}$. Then, $\overline{\varphi}_0$ is clearly $F$-invariant. Note that $\overline{\varphi}_0 = \varphi$ on $C_{bu}(\mathbb{R})$, which implies that $\overline{\varphi}_0 = \varphi$ from the $F$-invariance of $\overline{\varphi}_{0}$ and $\varphi$. We also note that $\|\overline{\varphi}_0\| = \|\varphi_0\|$ can be proved as above. Thus, we obtain
\[\|\varphi\| = \|\overline{\varphi}_0\| = \|\varphi_0\| = \|\varphi_{0, +}\| + \|\varphi_{0, -}\| = \|\overline{\varphi}_{0, +}\| + \|\overline{\varphi}_{0, -}\|.  \]
Therefore, by setting $\varphi_+ := \overline{\varphi}_{0, +}$ and $\varphi_- := \overline{\varphi}_{0, -}$, we obtain the desired decomposition of $\varphi$.
\end{prf}

\section{Iteration of sublinear functional $\overline{F}$}
In this section, we deal with sublinear functionals $\overline{F}_{\infty}$ induced by a functional $\overline{F}$ through infinite iteration. We need the following result from operator theory (see \cite{Kat}): for any contraction $U$ on a Banach space $X$, we put
\[\Gamma(U) = \sigma(U) \cap \Gamma, \]
where $\sigma(U)$ denotes the spectrum of $U$ and $\Gamma$ denotes the unit circle $\{z \in \mathbb{C} : |z| = 1\}$. Here, $\Gamma(U)$ is called the peripheral spectrum of $U$.
\begin{thm}
Let $U$ be a linear contraction on a Banach space $X$. Then, $\lim_{n \to \infty} \|U^n - U^{n+1}\| = 0$ if and only if the peripheral spectrum $\Gamma(U)$ of $U$ consists of at most the point $z =1$.
\end{thm}

By using this theorem and the results in Section 3, we can obtain the following result.
\begin{thm}
Let $f \in \mathscr{L}^1_{\flat}(\mathbb{R})$ and $F$ be the induced convolution operator. Then, 
\[\overline{F}_{\infty}(\phi) = \overline{P}(\phi) \]
holds for every $\phi \in L^{\infty}(\mathbb{R})$.
\end{thm}

\begin{prf}
First, note that the convolution operator $F$ induced by $f \in \mathscr{L}^1_{\flat}(\mathbb{R})$ satisfies $\|F\| = 1$ from the assumption that $\hat{f}(0) = 1$ and $f \ge 0$. Thus, $F$ is a contraction on $L^{\infty}(\mathbb{R})$. Furthermore, from Theorems 2.6 and 3.3, $\Gamma(F) = {1}$ holds. Thus, consideration of Theorem 4.1 yields 
\[\lim_{n \to \infty} \|F^n - F^{n+1}\| = 0. \]
Thus, for any $\phi \in L^{\infty}(\mathbb{R})$, we have
\[\lim_{n \to \infty} \|F^n\phi - F^{n+1}\phi\|_{\infty} = 0. \]
Additionally, for each $\phi \in L^{\infty}(\mathbb{R})$, 
\begin{align}
\overline{F}_{\infty}(\phi-F\phi) &= \lim_{n \to \infty} \overline{F}_n(\phi-F\phi) = \lim_{n \to \infty} \limsup_{x \to \infty} F^n(\phi(x)-F\phi(x)) \notag \\
&= \lim_{n \to \infty} \limsup_{x \to \infty} (F^n\phi(x) - F^{n+1}\phi(x)) \notag \\
&\le \lim_{n \to \infty} \|F^n\phi - F^{n+1}\phi\|_{\infty} = 0. \notag
\end{align}
\end{prf}
Note that, in the same way, $\underline{F}_{\infty}(\phi - f * \phi) = 0$ can also be proved. 
Let $\mathscr{P}$ and $\mathscr{F}_{\infty}$ be the weak*-closed convex subsets of $L^{\infty}(\mathbb{R})^*$; they are defined by $\varphi \in \mathscr{P}$ if and only if $\varphi(\phi) \le \overline{P}(\phi)$ for all $\phi \in L^{\infty}(\mathbb{R})$ and $\varphi \in \mathscr{F}_{\infty}$ if and only if $\varphi(\phi) \le \overline{F}_{\infty}(\phi)$ for all $\phi \in L^{\infty}(\mathbb{R})$, respectively. Since 
\[\underline{F}_{\infty}(\phi - f * \phi) \le \varphi(\phi - f * \phi) \le \overline{F}_{\infty}(\phi - f * \phi) \]
holds for each $\varphi \in \mathscr{F}_{\infty}$, for any $\varphi \in \mathscr{F}_{\infty}$, $\varphi$ is $F$-invariant. Therefore, from Theorem 3.7, $\varphi \in \mathscr{P}$ holds.
Hence, we have
\[\overline{F}_{\infty}(\phi) = \sup_{\varphi \in \mathscr{F}_{\infty}} \varphi(\phi) \le \sup_{\varphi \in \mathscr{P}} \varphi(\phi) = \overline{P}(\phi) \]
for every $\phi \in L^{\infty}(\mathbb{R})$. Conversely, from Corollary 3.2, note that the sublinear functional $\overline{P}(\phi)$ is $F$-invariant, i.e., $\overline{P}(\phi) = \overline{P}(F\phi) = \overline{P}(F^2\phi) = \cdots =\overline{P}(F^k\phi) = \cdots$ holds for every $\phi \in L^{\infty}(\mathbb{R})$. Hence, for each $k \ge 1$, we have
\[\overline{P}(\phi) = \overline{P}(F^k\phi) \le \limsup_{x \to \infty} F^k\phi(x) = \overline{F}_k(\phi). \]
Thus, we obtain
\[\overline{P}(\phi) \le \lim_{k \to \infty} \overline{F}_k(\phi) = \overline{F}_{\infty}(\phi). \]
Hence, we have obtained the desired equation
\[\overline{F}_{\infty}(\phi) = \overline{P}(\phi) \]
for each $\phi \in L^{\infty}(\mathbb{R})$.

In conjunction with Theorems 3.7 and 3.8, we obtain the following.
\begin{thm}
Let $f \in \mathscr{L}^1_{\flat}(\mathbb{R})$ and $F$ be the induced convolution operator. Then, for $\varphi \in L^{\infty}(\mathbb{R})^*$, $\varphi \in \mathscr{M}_F$ if and only if 
\[\varphi(\phi) \le \overline{P}(\phi) \]
holds for every $\phi \in L^{\infty}(\mathbb{R})$. For $\varphi \in M_F$, there exists unique elements $\varphi_+$ and $\varphi_-$ in $\mathscr{M}_F$ and nonnegative numbers $\alpha$ and $\beta$ such that
\[\varphi = \alpha \varphi_+ - \beta \varphi_-, \quad \|\varphi\| = \alpha\|\varphi_+\| + \beta\|\varphi_-\|. \]
Furthermore, for every $\phi \in L^{\infty}(\mathbb{R})$, the formula
\[\overline{F}_{\infty}(\phi) = \overline{P}(\phi) \]
holds.
\end{thm}

\section{Relationship between Wiener's Tauberian theorem and almost convergence}
In this section, we deal with the relationship between two summability methods on $L^{\infty}(\mathbb{R})$, i.e., summability methods defined via a convolution operator $F$ and a continuous analogue of almost convergence. For any $f \in L^1(\mathbb{R})$, let us define a summability method $F$ by
\[F(\phi) = \lim_{x \to \infty} \int_{-\infty}^{\infty} \phi(x-t)f(t)dt, \]
provided that the limit exists. Note that the functional $\underline{F}(\phi) := -\overline{F}(-\phi)$ on $L^{\infty}(\mathbb{R})$ can be expressed by
\[\underline{F}(\phi) = \liminf_{x \to \infty} \int_{-\infty}^{\infty} \phi(x-t)f(t)dt. \]
It is then obvious that $F(\phi) = \alpha$ if and only if $\overline{F}(\phi) = \underline{F}(\phi) = \alpha$. 
A particularly important case is when $f \in L^1(\mathbb{R})$ is a Wiener kernel, i.e., the zero set $Z(f)$ of the Fourier transform of $f$ is empty. The importance of Wiener kernels is illustrated by Wiener's Tauberian theorem, which is formulated as follows.
\begin{thm}
Let $f \in L^1(\mathbb{R})$ be a Wiener kernel with $\hat{f}(0) =1$. Let $g \in L^1(\mathbb{R})$ with $\hat{g}(0) = 1$. Then, for any $\phi \in L^{\infty}(\mathbb{R})$, if 
\[\lim_{x \to \infty} \int_{-\infty}^{\infty} \phi(x-t)f(t)dt = \alpha, \]
then 
\[\lim_{x \to \infty} \int_{-\infty}^{\infty} \phi(x-t)g(t)dt = \alpha. \]
\end{thm}
We now refer to the following simple observation that yields a simpler formulation of Wiener's Tauberian theorem.
\begin{thm}
Let $f \in L^1(\mathbb{R})$ be a Wiener kernel with $\hat{f}(0) =1$. Then, for any $\phi \in L^{\infty}(\mathbb{R})$, if 
\[\lim_{x \to \infty} \int_{-\infty}^{\infty} \phi(x-t)f(t)dt = \alpha, \]
then
\[w^*\mathchar`-\lim_{s} \phi_s(x) = \alpha, \]
where the symbol $w^*\mathchar`-\lim$ denotes the limit in the weak*-topology of $L^{\infty}(\mathbb{R})$.
\end{thm}

\begin{prf}
Since a limit is translation invariant, for any $s \in \mathbb{R}$, we have
\begin{align}
\lim_{x \to \infty} \int_{-\infty}^{\infty} \phi(x-t)f(t)dt = 0 &\Longleftrightarrow \lim_{x \to \infty} \int_{-\infty}^{\infty} \phi_x(-t)f(t)dt = 0 \notag \\
&\Longleftrightarrow \lim_{x \to \infty} \int_{-\infty}^{\infty} \phi_{x}(t+s)f(-t)dt = 0 \notag \\
&\Longleftrightarrow \lim_{x \to \infty} \int_{-\infty}^{\infty} \phi_x(t)f(s-t)dt = 0. \notag
\end{align}
Hence, for any element $h$ in the closed linear hull of the translates $\{f_s(-t)\}_{s \in \mathbb{R}}$, we have 
\[\lim_{x \to \infty} \int_{-\infty}^{\infty} \phi_x(t)h(t)dt = 0. \]
From the assumption that $f$ is a Wiener kernel and fromTheorem 2.4, the functions $h$ consist of all functions in $L^1(\mathbb{R})$, and the result follows immediately.
\end{prf}
This theorem implies that any summability method $F$ induced by a Wiener kernel $f \in L^1(\mathbb{R})$ with $\hat{f}(0) = 1$ is equivalent to the summability method $W$ defined by $W(\phi) = \alpha$ if and only if $w^*\mathchar`-\lim_{s} \phi_s(x) = \alpha$. Since the expression of $W$ has the advantage of being independent of specific Wiener kernels, we use it occasionally thereafter.

Now, we define a summability method via the functional $\overline{P}$. Note that the functional $\underline{P}(\phi) := -\overline{P}(-\phi)$ can be expressed by
\[\underline{P}(\phi) = \lim_{\theta \to \infty} \liminf_{x \to \infty} \frac{1}{\theta} \int_x^{x+\theta} \phi(t)dt. \]
Then, for $\phi \in L^{\infty}(\mathbb{R})$, we define $P(\phi) = \alpha$ if and only if $\overline{P}(\phi) = \underline{P}(\phi) = \alpha$. As the following theorems illustrate, this can be viewed as a continuous version of the classical notion of almost convergence introduced by Lorentz (\cite{Lor}).
\begin{thm}
Let $\phi \in L^{\infty}(\mathbb{R})$. Then, $P(\phi) = \alpha$ if and only if there exists a constant $R \ge 0$ for every $\varepsilon > 0$ such that if $\theta \ge R$, then
\[\left|\frac{1}{\theta} \int_x^{x+\theta} \phi(t)dt - \alpha \right| \le \varepsilon \]
for a sufficiently large $x \ge 0$.
\end{thm}

\begin{prf}
First, we prove the sufficiency. Suppose that the above assertion holds. We show that $\overline{P}(\phi) = \underline{P}(\phi) = \alpha$. For any fixed $\varepsilon > 0$, there exists $R \ge 0$ such that if $\theta \ge R$, then
\[\alpha - \varepsilon \le \frac{1}{\theta} \int_x^{x+\theta} \phi(t)dt \le \alpha + \varepsilon \]
for a sufficiently large $x \ge 0$. This implies that
\[\alpha - \varepsilon \le \liminf_{x \to \infty} \frac{1}{\theta} \int_x^{x+\theta} \phi(t)dt \le \limsup_{x \to \infty} \int_x^{x+\theta} \phi(t)dt \le \alpha + \varepsilon \]
whenever $\theta \ge R$. Since $\varepsilon > 0$ is arbitrary, we have
\[\lim_{\theta \to \infty} \liminf_{x \to \infty} \frac{1}{\theta} \int_x^{x+\theta} \phi(t)dt = \lim_{\theta \to \infty} \limsup_{x \to \infty} \int_x^{x+\theta} \phi(t)dt = \alpha. \]
The desired result is thus obtained. 

Next, we prove the necessity. Suppose that $P(\phi) = \alpha$, i.e., $\overline{P}(\phi) = \underline{P}(\phi) = \alpha$. Then, for any $\varepsilon > 0$, there exists a constant $R \ge 0$ such that 
\[\liminf_{x \to \infty} \frac{1}{\theta} \int_x^{x + \theta} \phi(t)dt \ge \alpha - \frac{\varepsilon}{2} \]
and 
\[\limsup_{x \to \infty} \frac{1}{\theta} \int_x^{x + \theta} \phi(t)dt \le \alpha + \frac{\varepsilon}{2} \]
whenever $\theta \ge R$. Furthermore, we can choose a constant $R_{\theta} \ge 0$ such that 
\[\alpha - \varepsilon \le \frac{1}{\theta} \int_x^{x+\theta} \phi(t)dt \le \alpha + \varepsilon \]
whenever $x \ge R_{\theta}$. Hence, for any $\varepsilon > 0$, there exists $R \ge 0$ such that if $\theta \ge R$, then
\[\left|\frac{1}{\theta} \int_x^{x+\theta} \phi(t)dt - \alpha \right| \le \varepsilon \]
for $x \ge R_{\theta}$. The proof is now complete.
\end{prf}
We note that this assertion is equivalent to the following apparently stronger condition.
\begin{thm}
Let $\phi \in L^{\infty}(\mathbb{R})$. Then, $P(\phi) = \alpha$ if and only if 
\[\lim_{\theta \to \infty} \frac{1}{\theta} \int_x^{x+\theta} \phi(t)dt = \alpha \]
holds uniformly in $x \ge 0$.
\end{thm}
The proof is the same as the proof of Theorem 14 in \cite{Suk}, where the assertion is proved for the elements of $C_{bu}(\mathbb{R})$. 

Now, we direct attention to the relationship between the summability methods $F$ and $P$. First, we show an Abelian theorem. The following result shows that almost convergence $P$ is stronger than any summability method $F$ induced by $f \in L^1(\mathbb{R})$ with $\hat{f}(0) = 1$.
\begin{thm}
Let $\phi \in L^{\infty}(\mathbb{R})$. Further, let $f \in L^1(\mathbb{R})$ satisfy $\hat{f}(0) = 1$. If
\[\lim_{x \to \infty} \int_{-\infty}^{\infty} \phi(t)f(x-t)dt = \alpha, \]
then
\[P(\phi) = \alpha. \]
\end{thm}
\begin{prf}
Let $\phi \in L^{\infty}(\mathbb{R})$. From Corollary 3.2, we have
\[\overline{P}(\phi) = \overline{P}(f * \phi) \le \limsup_{x \to \infty} (f * \phi)(x) = \overline{F}(\phi). \]
Further, we have
\[\underline{P}(\phi) = -\overline{P}(-\phi) = -\overline{P}(-f*\phi) \ge \liminf_{x \to \infty} (f * \phi)(x) = \underline{F}(\phi). \]
Thus, we obtain the following relation:
\[\underline{F}(\phi) \le \underline{P}(\phi) \le \overline{P}(\phi) \le \overline{F}(\phi). \]
This implies the required result.
\end{prf}

Next, we consider a Tauberian theorem, i.e., we show that under a certain condition, $P(\phi) = \alpha$ implies $F(\phi) = \alpha$. We begin with the following lemma.
\begin{lem}
Let $\phi \in L^{\infty}(\mathbb{R})$. Then, the following conditions are equivalent:
\begin{enumerate}
\item For every $s \in \mathbb{R}$, $w^*\mathchar`-\lim_{x} (\phi_{x+s} - \phi_x) = 0$, i.e., $W(\phi_s -\phi) = 0$;
\item $\int_{-\infty}^{\infty} \phi(x-t)f(t)dt = 0$ for an $f \in L^1(\mathbb{R})$ with $Z(f) = \{0\}$.
\end{enumerate}
\end{lem}

\begin{prf}
(1) $\Rightarrow$ (2): Let $f \in L^1(\mathbb{R})$ be such that $Z(f) = \emptyset$, and let $s$ be a real number. Then, from assumption (1), we have
\begin{align}
&\lim_{x \to \infty} \int_{-\infty}^{\infty} \{\phi_{x+s}(t) - \phi_x(t)\}f(-t)dt = 0 \notag \\
&\Longleftrightarrow \lim_{x \to \infty} \int_{-\infty}^{\infty} \phi_x(t)\{f(s-t) - f(-t)\}dt = 0 \notag \\
&\Longleftrightarrow \lim_{x \to \infty} \int_{-\infty}^{\infty} \phi(x-t)\{f_s(t) - f(t)\}dt = 0. \notag
\end{align}
Substitute $g_s = f_s - f$ and then $\hat{g_s}(\xi) = (e^{is\xi} -1)\hat{f}(\xi)$. Since we have assumed that $Z(f) = \emptyset$ and that $s$ is arbitrary, the zero set of the closed linear hull of the translates of functions $\{g_s\}_{s \in \mathbb{R}}$ is $\{0\}$. Since $\{0\}$ is a spectral synthesis set, for any $h \in L^1(\mathbb{R})$ with $Z(h) = \{0\}$, we have $H(\phi) = 0$; thus, (2) holds.

(2) $\Rightarrow$ (1): This can be proved in the same way as above.
\end{prf}

\begin{thm}
Let $\phi \in L^{\infty}(\mathbb{R})$. Then $W(\phi) = \alpha$ if and only if $P(\phi) = \alpha$ and one of the two conditions in Lemma 5.1 holds. In other words, either condition of Lemma 5.1 is a Tauberian condition under which $P$ summability implies $F$ summability.
\end{thm}

\begin{prf}
The necessity is self-evident. Thus, we prove the sufficiency.
Note that from Theorem 5.1, without loss of generality, we can assume that the Wiener kernel $f$ is in $\mathscr{L}^1_{\flat}(\mathbb{R})$. Observe that
\[\int_{-\infty}^{\infty} \{\phi(t) - (F\phi)(t)\}f(x-t)dt = \int_{-\infty}^{\infty} \phi(t)\{f(x-t) - f^{*2}(x-t)\}dt. \]
Substitute $g = f - f^{*2}$. Then, we have $\hat{g}(\xi) = \hat{f}(\xi) (1 - \hat{f}(\xi))$. Thus, we obtain $Z(g) = \{0\}$. From the assumption, the right-hand side of the above equation is $0$. Hence, we have $F(\phi - F\phi) = 0$, which implies that $F(F^k\phi - F^{k+1}\phi) = 0$ for every $k \ge 1$. From Theorem 4.2, we have
\begin{align}
\overline{F}(\phi) &= \overline{F}(\phi - F\phi + F\phi - F^2\phi + F^2\phi - \cdots - F^k\phi + F^k\phi)  \notag \\
&= F(\phi - F\phi) + F(F\phi - F^2\phi) + \cdots F(F^{k-1}\phi - F^k\phi) + \overline{F}(F^k\phi) \notag \\ 
&= \overline{F}(F^k\phi) = \overline{F}_{k+1}(\phi) \notag
\end{align}
for every $k \ge 1$. Hence, we obtain
\[\overline{F}(\phi) = \lim_{k \to \infty} \overline{F}_k(\phi) = \overline{P}(\phi) = \alpha. \]
In the same way, we have
\[\underline{F}(\phi) = \lim_{k \to \infty} \underline{F}_k(\phi) = \underline{P}(\phi) = \alpha. \]
Therefore, we obtain the result $F(\phi) = \overline{F}(\phi) = \underline{F}(\phi) = \alpha$.
\end{prf}

\section{Mellin convolution-invariant functionals on $L^{\infty}(\mathbb{R}^{\times})$}
In this section, we deal with a multiplicative version of the results in the preceding sections. Let $L^{\infty}(\mathbb{R}^{\times})$ be the set of all essentially bounded functions on $\mathbb{R}^{\times}$ and let $L^1(\mathbb{R}^{\times})$ be the group algebra of $\mathbb{R}^{\times}$. Note that the Haar measure of the positive multiplicative group $\mathbb{R}^{\times}$ is $\frac{dt}{t}$. For any $g \in L^1(\mathbb{R}^{\times})$, let the symbol $G$ denote the convolution operator as follows:
\[G : L^{\infty}(\mathbb{R}^{\times}) \rightarrow L^{\infty}(\mathbb{R}^{\times}), \quad (G\phi)(x) = (g \overset{M}{*} \phi)(x), \quad \phi \in L^{\infty}(\mathbb{R}^{\times}), \]
where the Mellin convolution $\overset{M}{*}$ of $\phi \in L^{\infty}(\mathbb{R}^{\times})$ and $g \in L^1(\mathbb{R}^{\times})$ is defined by
\[(g \overset{M}{*} \phi)(x) = \int_0^{\infty} \phi(x/t)g(t)\frac{dt}{t} = \int_0^{\infty} \phi(t)g(x/t)\frac{dt}{t}. \]
Let $\mathscr{L}^1_*(\mathbb{R}^{\times})$ be the set of functions $g$ in $L^1(\mathbb{R}^{\times})$ such that $\hat{g}(\xi) = 1$ only at the point $\xi = 0$. Here, $\hat{g}$ is the Fourier transform of $g$ and is defined by
\[\hat{g}(x) = \int_0^{\infty} g(t)t^{ix} \frac{dt}{t}, \]
where $t^{ix} = e^{ix\log t}$. Let $\mathscr{L}^1_{\flat}(\mathbb{R}^{\times})$ be the set of functions $g$ in $L^1(\mathbb{R}^{\times})$ such that $g \ge 0$ and $\hat{g}(0)=1$. Further, we define $\mathscr{L}^1_{\sharp}(\mathbb{R}^{\times})$ as the set of functions $g$ in $L^1(\mathbb{R}^{\times})$ such that $g \ge 0$, $\hat{g}(0) = 1$, and $|\hat{g}(\xi)| = 1$ only at the point $\xi = 0$. The relation $\mathscr{L}^1_{\sharp}(\mathbb{R}^{\times}) \subseteq \mathscr{L}^1_*(\mathbb{R}^{\times})$ follows immediately from the definitions. Further, we have $\mathscr{L}^1_{\flat}(\mathbb{R}^{\times}) = \mathscr{L}^1_{\sharp}(\mathbb{R}^{\times})$, the proof of which is the same as that of Theorem 2.6. Thus, we can also show that $g \in \mathscr{L}^1_{\flat}(\mathbb{R}^{\times})$ if and only if $g \in \mathscr{L}^1_*(\mathbb{R}^{\times})$ and $g \ge 0$.

Let $L^{\infty}(\mathbb{R}^{\times})^*$ be the dual space of $L^{\infty}(\mathbb{R}^{\times})$ and $G^*$ be the adjoint operator of $G$:
\[G^* : L^{\infty}(\mathbb{R}^{\times})^* \rightarrow L^{\infty}(\mathbb{R}^{\times})^*, \quad (G^*\varphi)(\phi) = \varphi(G\phi). \]

Now, we consider $G$-invariant linear functionals, i.e., $\varphi \in L^{\infty}(\mathbb{R}^{\times})^*$ with $G^*\varphi = \varphi$, which vanish on $L^{\infty}(\mathbb{R}^{\times})_{0, +}$. Let us denote by $M_G$ the set of all $G$-invariant functionals and by $\mathscr{M}_G$ the set of all $G$-invariant means, which is the subset of $M_G$ with elements that satisfy the conditions $\varphi \ge 0$ and $\|\varphi\| =1$. 

An important fact is that through the results on the additive group of $\mathbb{R}$, we can obtain analogous results for Sections 3, 4, and 5 in the multiplicative setting.

Let us define the isometry $W$ of $L^{\infty}(\mathbb{R}^{\times})$ onto $L^{\infty}(\mathbb{R})$ as follows:
\[W : L^{\infty}(\mathbb{R}^{\times}) \rightarrow L^{\infty}(\mathbb{R}), \quad (W\phi)(x) = \phi(e^x). \]
Let $W^* : L^{\infty}(\mathbb{R})^* \rightarrow L^{\infty}(\mathbb{R}^{\times})^*$ be its adjoint operator. For any $f \in L^1(\mathbb{R})$, let us define a function $g$ in $L^1(\mathbb{R}^{\times})$ by  $g(x) = f(\log x)$. 
The following commutative diagram is then obtained and can be proved easily through integration by substitution:
$$
\begin{CD}
L^{\infty}(\mathbb{R}^{\times}) @>G >> L^{\infty}(\mathbb{R}^{\times})\\
@VWVV @VVWV \\
L^{\infty}(\mathbb{R}) @>F >> L^{\infty}(\mathbb{R})
\end{CD}
$$
Notice that the mapping $L^1(\mathbb{R}) \ni f(x) \mapsto g = f(\log x) \in L^1(\mathbb{R}^{\times})$ preserves the order and spectral properties, that is, $f  \in \mathscr{L}^1_*(\mathbb{R})$ if and only if $g \in \mathscr{L}^1_*(\mathbb{R}^{\times})$, $f \in \mathscr{L}^1_{\sharp}(\mathbb{R})$ if and only if $g \in \mathscr{L}^1_{\sharp}(\mathbb{R}^{\times})$, and $f \in \mathscr{L}^1_{\flat}(\mathbb{R})$ if and only if $g \in \mathscr{L}^1_{\flat}(\mathbb{R}^{\times})$.

Assume that $\varphi \in L^{\infty}(\mathbb{R})^*$ is $F$-invariant. Then, $W^*\varphi$ is $G$-invariant. In fact, since $G = W^{-1}FW$ from the above diagram, we have
\[(G^*(W^*\varphi))(\phi) = (W^*\varphi)(G\phi) = \varphi(WG\phi) = \varphi(FW\phi) = \varphi(W\phi) = (W^*\varphi)(\phi) \]
for every $\phi \in L^{\infty}(\mathbb{R}^{\times})$. Thus, $G^*W^*\varphi = W^*\varphi$, and $W^*\varphi$ is $G$-invariant. Consequently, we have proved the following result: the restriction of $W^*$ to $M_F$,
\[W^* : M_F \rightarrow M_G, \]
is a linear isomorphism between $M_F$ and $M_G$. It is also clear that the restriction of $W^*$ to $\mathscr{M}_F$,
\[W^* : \mathscr{M}_F \rightarrow \mathscr{M}_G, \]
is an affine isomorphism between $\mathscr{M}_F$ and $\mathscr{M}_G$.

Let $\overline{Q}: L^{\infty}(\mathbb{R}^{\times}) \rightarrow \mathbb{R}$ be the sublinear functional defined by
\[\overline{Q}(\phi) = \lim_{\theta \to \infty} \limsup_{x \to \infty} \frac{1}{\log \theta} \int_x^{\theta x} \phi(t)\frac{dt}{t}, \quad \phi \in L^{\infty}(\mathbb{R}^{\times}). \]
We can then easily confirm the following result through direct computation (integration by substitution).
\begin{lem}
For every $\phi \in L^{\infty}(\mathbb{R}^{\times})$, 
\[\overline{P}(W\phi) = \overline{Q}(\phi) \] 
holds.
\end{lem}
We provide a multiplicative version of Theorems 3.7, 3.8, and 4.2 as follows:

\begin{thm}
Let $g \in \mathscr{L}^1_*(\mathbb{R}^{\times})$ and $G$ be the induced operator. Then, for $\varphi \in L^{\infty}(\mathbb{R}^{\times})^*$, $\varphi \in \mathscr{M}_G$ if and only if
\[\varphi(\phi) \le \overline{Q}(\phi) \]
holds for every $\phi \in L^{\infty}(\mathbb{R}^{\times})$.
\end{thm}

\begin{prf}
It is sufficient to prove that $\sup_{\psi \in \mathscr{M}_G} \psi(\phi) = \overline{Q}(\phi)$ for every $\phi \in L^{\infty}(\mathbb{R}^{\times})$.
Note that for any $\psi \in \mathscr{M}_G$, there exists a $\varphi \in \mathscr{M}_F$ such that $W^*\varphi = \psi$ from the fact mentioned above. Thus, from Theorem 3.7 and Lemma 6.1, for any $\phi \in L^{\infty}(\mathbb{R}^{\times})$, we have 
\[\sup_{\psi \in \mathscr{M}_G} \psi(\phi) = \sup_{\varphi \in \mathscr{M}_F} (W^*\varphi)(\phi) = \sup_{\varphi \in \mathscr{M}_F} \varphi(W\phi) = \overline{P}(W\phi) = \overline{Q}(\phi). \]
\end{prf}

The following is a corresponding expression for Theorem 3.8 that can be easily deduced via the isomorphism $W^*$ of $M_F$ onto $M_G$.
\begin{thm}
Let $G$ be a convolution operator induced by an element $g$ of $\mathscr{L}^1_{\flat}(\mathbb{R}^{\times})$. Then, for any $\varphi \in M_G$, $\varphi$ admits the Jordan decomposition in $M_G$, i.e., there exist positive elements $\varphi_+$ and $\varphi_-$ in $M_G$ such that 
\[\varphi = \varphi_+ - \varphi_- , \quad \|\varphi\| = \|\varphi_+\| + \|\varphi_-\|.    \]
\end{thm}

The corresponding result for Theorem 4.2 is as follows.
\begin{thm}
Let $g \in \mathscr{L}^1_{\flat}(\mathbb{R}^{\times})$ and $G$ be the induced convolution operator . Then, 
\[\overline{G}_{\infty}(\phi) = \overline{Q}(\phi) \]
holds for every $\phi \in L^{\infty}(\mathbb{R}^{\times})$.
\end{thm}

\begin{prf}
We only show that $\overline{G}_{\infty}(\phi-G\phi) = 0$ for every $\phi \in L^{\infty}(\mathbb{R}^{\times})$. The remainder of the proof is similar to that for Theorem 4.2 and is left to the reader. Observe that
\begin{align}
\overline{G}_{\infty}(\phi - G\phi) &= \lim_{n \to \infty} \overline{G}_n(\phi-G\phi) = \lim_{n \to \infty} \limsup_{x \to \infty} (G^n\phi(x) - G^{n+1}\phi(x)) \notag \\
&= \lim_{n \to \infty} \limsup_{x \to \infty} ((W^{-1}F^nW\phi)(x) - (W^{-1}F^{n+1}W\phi)(x)) \notag \\
&= \lim_{n \to \infty} \limsup_{x \to \infty} ((F^nW\phi)(\log x) - (F^{n+1}W\phi)(\log x)) \notag \\
&\le \lim_{n \to \infty} \|F^nW\phi - F^{n+1}W\phi\|_{\infty} = 0. \notag
\end{align}
Note that in the last equation, we use a result from the proof of Theorem 4.2.
\end{prf}

In conjunction with Theorems 6.1 and 6.2, we have the following.
\begin{thm}
Let $g \in \mathscr{L}^1_{\flat}(\mathbb{R}^{\times})$ and $G$ be the induced convolution operator. Then, for $\varphi \in L^{\infty}(\mathbb{R}^{\times})^*$, $\varphi \in \mathscr{M}_G$ if and only if 
\[\varphi(\phi) \le \overline{Q}(\phi) \]
holds for every $\phi \in L^{\infty}(\mathbb{R}^{\infty})$. For $\varphi \in M_G$, there exists unique elements $\varphi_+$ and $\varphi_-$ in $\mathcal{M}_G$ and nonnegative numbers $\alpha$ and $\beta$ such that
\[\varphi = \alpha \varphi_+ - \beta \varphi_-, \quad \|\varphi\| = \alpha\|\varphi_+\| + \beta\|\varphi_-\|. \]
Further, for every $\phi \in L^{\infty}(\mathbb{R}^{\times})$, the equation
\[\overline{G}_{\infty}(\phi) = \overline{Q}(\phi) \]
holds.
\end{thm}

Now, let us consider an example. For $r > 0$, let us define $g_r(x) \in \mathscr{L}^1_{\flat}(\mathbb{R}^{\times})$ by 
\[g_r(x) = \begin{cases}
            rx^{-r} & \text{if $x \ge 1$}, \\
            0       & \text{if $x < 1$}.
         \end{cases}     
\]
The corresponding convolution operator $G_r : L^{\infty}(\mathbb{R}^{\times}) \rightarrow L^{\infty}(\mathbb{R}^{\times})$ is given as follows:
\[(G_r\phi)(x) = \int_0^{\infty} \phi(t) r\left(\frac{x}{t}\right)^{-r}\frac{dt}{t} = \frac{r}{x^r} \int_0^x \phi(t)t^{r-1}dt, \]
where $\phi \in L^{\infty}(\mathbb{R}^{\times})$. In particular, for $r=1$, we have the Hardy operator $G_1 = H$.

From Theorem 6.4, we obtain the following result, which includes the main theorems of \cite{Kuni} as a special case.
\begin{thm}
For $\varphi \in L^{\infty}(\mathbb{R}^{\times})^*$, $\varphi$ is a $G_r$-invariant mean if and only if 
\[\varphi(\phi) \le \overline{Q}(\phi) \]
holds for every $\phi \in L^{\infty}(\mathbb{R}^{\times})$. Further, 
\[\overline{G}_{r, \infty}(\phi) = \overline{Q}(\phi) \]
for every $\phi \in L^{\infty}(\mathbb{R}^{\times})$.
\end{thm}

Finally, we present an application to summability methods. In what follows, we present only assertions without proofs since they are similar to those in Section 5. For any $g \in L^1(\mathbb{R}^{\times})$, let us define the summability method $G$ by
\[G(\phi) = \lim_{x \to \infty} \int_0^{\infty} \phi(x/t)g(t) \frac{dt}{t}, \]
provided that the limit exists. For $\phi \in L^{\infty}(\mathbb{R}^{\times})$ and $r > 0$, let us define $\phi^*_r(x) = \phi(rx)$. The multiplicative analogue of Theorem 5.2 is as follows.
\begin{thm}
Let $g \in L^1(\mathbb{R}^{\times})$ satisfy $Z(g) = \emptyset$ and $\hat{g}(0) = 1$. Then, if
\[\lim_{x \to \infty} \int_0^{\infty} \phi(x/t)g(t) \frac{dt}{t} = \alpha \]
for any $\phi \in L^{\infty}(\mathbb{R}^{\times})$, then
\[w^*\mathchar`-\lim_{r} \phi^*_r(x) = \alpha, \]
where the symbol $w^*\mathchar`-\lim$ denotes the limit in the weak*-topology of $L^{\infty}(\mathbb{R}^{\times})$.
\end{thm}
Let $W^*$ be the summability method on $L^{\infty}(\mathbb{R}^{\times})$ defined by $W^*(\phi) = \alpha$ if and only if $w^*\mathchar`-\lim_r \phi^*_r(x) = \alpha$. Further, let us define the summability method $Q$ on $L^{\infty}(\mathbb{R}^{\times})$ by $Q(\phi) = \alpha$ if and only if $\overline{Q}(\phi) = \underline{Q}(\phi) = \alpha$, where $\underline{Q}(\phi) = -\overline{Q}(-\phi)$. Then, the following result, which corresponds to Theorem 5.4, is obtained.
\begin{thm}
Let $\phi \in L^{\infty}(\mathbb{R}^{\times})$. Then, $Q(\phi) = \alpha$ if and only if
\[\lim_{\theta \to \infty} \frac{1}{\log \theta} \int_x^{\theta x} \phi(t) \frac{dt}{t} = \alpha \]
holds uniformly in $x \ge 1$.
\end{thm}

Subsequently, we address Abelian and Tauberian theorems between summability methods $G$ and $Q$. Note that the following results are counterparts of Theorem 5.5, Lemma 5.1, and Theorem 5.6, respectively.
\begin{thm}
Let  $\phi \in L^{\infty}(\mathbb{R}^{\times})$. Further, let $g \in L^1(\mathbb{R}^{\times})$ satisfy $\hat{g}(0) = 1$. If
\[\lim_{x \to \infty} \int_0^{\infty} \phi(t)g(x/t) \frac{dt}{t} = \alpha, \]
then
\[Q(\phi) = \alpha. \]
\end{thm}

\begin{lem}
Let $\phi \in L^{\infty}(\mathbb{R}^{\times})$. Then, the following conditions are equivalent:
\begin{enumerate}
\item For every $r > 0$, $w^*\mathchar`-\lim_{x} (\phi^*_{rx} - \phi^*_x) = 0$, i.e., $W^*(\phi^*_r - \phi) = 0$;
\item $\int_0^{\infty} \phi(x/t)g(t)dt = 0$ for a $g \in L^1(\mathbb{R}^{\times})$ with $Z(g) = \{0\}$. 
\end{enumerate}
\end{lem}

\begin{thm}
Let $\phi \in L^{\infty}(\mathbb{R}^{\times})$. Then, $W^*(\phi) = \alpha$ if and only if $Q(\phi) = \alpha$ and one of the two conditions in Lemma 6.2 holds. In other words, either condition of Lemma 6.2 is a Tauberian condition under which $Q$ summability implies $G$ summability.
\end{thm}

\section{Applications to Ces\`{a}ro operators}
In this section, we consider the means on $l_{\infty}$ of bounded functions on natural numbers $\mathbb{N}$ invariant with respect to the Ces\`{a}ro operator and their application to summability methods on functions on $\mathbb{N}$. Recall that the Ces\`{a}ro operator is defined as follows:
\[C : l_{\infty} \rightarrow l_{\infty}, \quad (C\phi)(n) = \frac{1}{n} \sum_{i=1}^{n} \phi(i), \]
where $\phi \in l_{\infty}$. Ces\`{a}ro invariant means are defined as the elements $\varphi$ of $l_{\infty}^*$ of the dual space of $l_{\infty}$ that satisfy $\varphi \ge 0$, $\|\varphi\|=1$, and $C^*\varphi = \varphi$, where $C^*$ is the adjoint operator of $C$. It is obvious that the Ces\`{a}ro operator can be viewed as a discrete version of the Hardy operator $H$ on $L^{\infty}(\mathbb{R}^{\times})$, which is defined as follows:
\[H : L^{\infty}(\mathbb{R}^{\times}) \rightarrow L^{\infty}(\mathbb{R}^{\times}), \quad (H\phi)(x) = \frac{1}{x} \int_0^x \phi(t)dt. \]

We now consider the relationship between Ces\`{a}ro-invariant functionals $M_C$ and Hardy-invariant functionals $M_H$. Specifically, we show that $M_C$ and $M_H$ are isomorphic, and based on this result, we obtain results similar to Theorem 6.5 pertaining to the Ces\`{a}ro-invariant means, which provide another proof of the results in \cite{Kuni, Sem}. We also present an application of these results to the $C_{\infty}$ summability method, which is a generalization of the Ces\`{a}ro or  H\"{o}lder summability methods.

For a function $\phi$ on $\mathbb{N}$, recall that its Ces\`{a}ro mean $C(\phi)$ is defined by the limit
\[C(\phi) := \lim_{n \to \infty} \frac{1}{n} \sum_{i=1}^n \phi(i), \]
provided that the limit exists. This is a summability method of the simplest type. One of the simplest generalizations of the Ces\`{a}ro mean is its iteration, i.e., by using  the Ces\`{a}ro operator, the summability methods $C_1, C_2, \ldots C_k, \ldots$ are defined by
\[C_1(\phi) := C(\phi), C_2(\phi) := C(C\phi), \ldots, C_k(\phi) := C_{k-1}(C\phi), \ldots. \]
These classical methods are called H\"{o}lder summability methods (see \cite{Hardy}).
Furthermore, we can define the $C_{\infty}$ summability method, which was originally introduced by \cite{Gar}, by the limit of the sequence $\{C_k\}_{k \ge 1}$ of the H\"{o}lder summability methods. Let us define the sublinear functionals $\{\overline{C}_k\}_{k \ge 1}$ on $l_{\infty}$ as
\[\overline{C}_1(\phi) := \overline{C}(\phi) := \limsup_{n \to \infty} \frac{1}{n} \sum_{i=1}^n \phi(i), \]
and for $k = 2, 3, \cdots$,
\[\overline{C}_k(\phi) := \overline{C}_{k-1}(C\phi) = \overline{C}(C^{k-1}\phi). \]
Further, we define the lower version $\underline{C}_k$ of $\overline{C}_k$ by
\[\underline{C}_1(\phi) := \underline{C}(\phi) := \liminf_{n \to \infty} \frac{1}{n} \sum_{i=1}^n \phi(i), \]
and for $k = 2, 3, \cdots$,  
\[\underline{C}_k(\phi) := \underline{C}_{k-1}(C\phi) = \underline{C}(C^{k-1}\phi). \]
It is obvious that for every $\phi \in l_{\infty}$,
\[\overline{C}(\phi) \ge \overline{C}_2(\phi) \ge \ldots \ge \overline{C}_n(\phi) \ge \ldots \ge \underline{C}_n(\phi) \ge \ldots \underline{C}_2(\phi) \ge \underline{C}(\phi). \]
Thus, let us define the functionals $\overline{C}_{\infty}$ and $\underline{C}_{\infty}$ by
\[\overline{C}_{\infty}(\phi) = \lim_{k \to \infty} \overline{C}_k(\phi), \quad \underline{C}_{\infty}(\phi) = \lim_{k \to \infty} \underline{C}_k(\phi), \]
respectively. The $C_{\infty}$ summability method is then defined as follows: for a function $\phi$ on $\mathbb{N}$, $C_{\infty}(\phi) = \alpha$ if and only if $\overline{C}_{\infty}(\phi) = \underline{C}_{\infty}(\phi) = \alpha$.

This summability method was studied by several researchers \cite{DiaD, Dia, Dur, Ebe, Fle, Gar}. \cite{Sem} also dealt with the sublinear functional $\overline{C}_{\infty}$, though not from the perspective of summability methods but rather from Ces\`{a}ro-invariant means. 

We obtain an analytic expression of the sublinear functional $\overline{C}_{\infty}$, from which we can deduce some of the properties of the $C_{\infty}$ summability method, including a necessary and sufficient condition that a given $\phi \in l_{\infty}$ is $C_{\infty}$ summable.

In what follows, the terms $H$-invariant and $C$-invariant mean Hardy invariant and Ces\`{a}ro invariant, respectively. We first show that $M_H$ and $M_C$ are isomorphic. Let us define the linear operator $V$ as follows: 
\[V : l_{\infty} \rightarrow L^{\infty}(\mathbb{R}^{\times}), \quad (V\phi)(x) = \phi([x+1]), \quad x > 0. \]
Let us consider its adjoint operator $V^* : L^{\infty}(\mathbb{R}^{\times})^* \rightarrow l_{\infty}^*$ and show that $V^*$ is a one-to-one and onto mapping from $M_H$ to $M_C$.

First, we show that if $\varphi \in L^{\infty}(\mathbb{R}^{\times})^*$ is $H$-invariant, then $V^*\varphi$ is $C$-invariant. Suppose that $\varphi \in M_H$. We show that $C^*V^*\varphi(\phi) = V^*\varphi(\phi)$ for every $\phi \in l_{\infty}$. Note that $C^*V^*\varphi(\phi) = V^*\varphi(C\phi) = \varphi(VC\phi)$, and
\begin{align}
(VC\phi)(x) &= (C\phi)([x+1]) = \frac{1}{[x+1]} \sum_{i=1}^{[x+1]} \phi(i) \notag \\
&= \frac{1}{[x+1]} \int_0^{[x+1]} (V\phi)(t)dt =  (HV\phi)([x+1]) \notag 
\end{align}
holds. Observe that 
\[\left|(HV\phi)(x) - (HV\phi)([x+1])\right| \rightarrow 0 \ \text{as} \ x \to \infty. \]
Then, from the assumption that $\varphi$ vanishes on $L^{\infty}_{0,+}(\mathbb{R})$, we have
\[C^*V^*\varphi(\phi) = \varphi(VC\phi) = \varphi(HV\phi) = \varphi(V\phi) = V^*\varphi(\phi). \]
Thus, we have proved that $V^*$ maps $M_H$ into $M_C$. 

Next, we show that $V^* : M_H \rightarrow M_C$ is surjective. Let us define the linear operator $V_1$ as
\[V_1 : L^{\infty}(\mathbb{R}^{\times}) \rightarrow l_{\infty}, \quad (V_1\phi)(n) = \int_{n-1}^n \phi(x)dx, \quad n \ge 1. \]
Let $V_1^* : l_{\infty}^* \rightarrow L^{\infty}(\mathbb{R}^{\times})^*$ be its adjoint operator. We show that if $\psi \in l_{\infty}^*$ is $C$-invariant, then $V_1^*\psi$ is $H$-invariant. Note that $H^*V_1^*\psi(\phi) = V_1^*\psi(H\phi) = \psi(V_1H\phi)$. First, observe that
\begin{align}
CV_1\phi(n) &= \frac{1}{n} \sum_{i=1}^n (V_1\phi)(i) = \frac{1}{n} \sum_{i=1}^n \int_{i-1}^i \phi(x)dx \notag \\
&= \frac{1}{n} \int_0^n \phi(t)dt = (H\phi)(n). \notag
\end{align}
Note that 
\begin{align}
|(V_1H\phi)(n) - (H\phi)(n)| &= \left|\int_{n-1}^n (H\phi)(t)dt - (H\phi)(n)\right| \notag \\
&= \left|\int_{n-1}^n \{(H\phi)(t) - (H\phi)(n)\}dt\right| \notag \\
&\le \sup_{x \in [n-1, n]}|(H\phi)(x)-(H\phi)(n)| \notag \\
&\le \sup_{x \in [n-1, n]} \frac{1}{n} |H\phi(x)| + \frac{1}{n}\|\phi\|_{\infty}, \notag
\end{align}
which tends to $0$ as $n \to \infty$. Thus, we have
\begin{align}
|(CV_1\phi)(n) - (V_1H\phi)(n)| &\le |(CV_1\phi)(n) - (H\phi)(n)| + |(H\phi)(n) - (V_1H\phi)(n)| \notag \\
&\le \sup_{x \in [n-1, n]} \frac{1}{n} |H\phi(x)| + \frac{1}{n}\|\phi\|_{\infty}. \notag
\end{align}
Since the last term tends to $0$ as $n \to \infty$, $CV_1\phi(n) - V_1H\phi(n)$ tends to $0$ as $n \to \infty$. Therefore, we have
\[H^*V_1^*\psi(\phi) = \psi(V_1H\phi) = \psi(CV_1\phi) = \psi(V_1\phi) = (V_1^*\psi)(\phi). \]
We have shown that $V_1^*\psi$ is $H$-invariant if $\psi$ is $C$-invariant. Moreover, notice that $V^*(V_1^*\psi) = \psi$ for every $\psi \in l_{\infty}^*$. In fact, for any $\phi \in l_{\infty}$, we have
\begin{align}
(V_1V\phi)(n) &= \int_{n-1}^{n} (V\phi)(t)dt = \int_{n-1}^n \phi([t+1])dt \notag \\
&= \int_{n-1}^n \phi(n)dt = \phi(n). \notag
\end{align}
This implies the required result immediately.

Now, given any $\psi \in M_C$, $\varphi = V_1^*\psi$ is in $M_H$ and $V^*\varphi = V^*V_1^*\psi = \psi$. Thus, $V^* : M_H \rightarrow M_C$ is surjective.

Finally, we show that $V^* : M_H \rightarrow M_C$ is injective. Suppose that $V^*\varphi = V^*\varphi_1$ for $\varphi, \varphi_1 \in M_H$. It is sufficient to show that for any $\phi \in L^{\infty}(\mathbb{R}^{\times})$, there exists a $\psi \in l_{\infty}$ such that $\lim_{x \to \infty} |(H\phi)(x) - (HV\psi)(x)| = 0$. In fact, if this holds, we have $\varphi(\phi) = \varphi(H\phi) = \varphi(HV\psi) = \varphi(V\psi) = V^*\varphi(\psi)$ and, similarly, $\varphi_1(\phi) = V^*\varphi_1(\psi)$. From the assumption, we have $\varphi(\phi) = \varphi_1(\phi)$ for any $\phi \in L^{\infty}(\mathbb{R}^{\times})$, and we obtain $\varphi = \varphi_1$.

If we set
\[\psi(n) := \int_{n-1}^n \phi(t)dt, \quad n \ge 1, \]
we have
\begin{align}
(HV\psi)(x) &= \frac{1}{x} \int_0^x (V\psi)(t)dt = \frac{1}{x} \sum_{i=1}^{[x]} \int_{i-1}^i (V\psi)(t)dt + \frac{1}{x} \int_{[x]}^x (V\psi)(t)dt \notag \\ 
&= \frac{1}{x} \sum_{i=1}^{[x]} \psi(i) + \frac{1}{x} \int_{[x]}^x \psi([x]+1)dt \notag \\
&= \frac{1}{x} \int_0^{[x]} \phi(t)dt + \frac{x - [x]}{x} \cdot \psi([x]+1) \notag \\
&= (H\phi)(x) - \frac{1}{x} \int_{[x]}^x \phi(t)dt + \frac{x - [x]}{x} \cdot \psi([x]+1). \notag
\end{align}
Since the second and third terms of the last expression tend to $0$ as $x \to \infty$, $\psi$ is the desired function. This completes the proof. 

In particular, the restriction of $V^*$ to $\mathscr{M}_H$ is an affine isomorphism between $\mathscr{M}_H$ and $\mathscr{M}_C$, where $\mathscr{M}_H$ and $\mathscr{M}_C$ are the sets of $H$-invariant and $H_C$-invariant means, respectively. Now, we consider the maximal value attained by $\mathscr{M}_C$ for any fixed $\phi \in l_{\infty}$:
\[\overline{p}_C(\phi) := \sup_{\psi \in \mathscr{M}_C} \psi(\phi). \]
According to the above observation and Theorem 6.5, we have
\begin{align}
\sup_{\psi \in \mathscr{M}_C} \psi(\phi) &= \sup_{\varphi \in \mathscr{M}_H} (V^*\varphi)(\phi) = \sup_{\varphi \in \mathscr{M}_H} \varphi(V\phi) \notag \\
&= \lim_{\theta \to \infty} \limsup_{x \to \infty} \frac{1}{\log \theta} \int_x^{\theta x} (V\phi)(t)\frac{dt}{t} \notag \\
&= \lim_{\theta \to \infty} \limsup_{n \to \infty} \frac{1}{\log \theta} \int_n^{\theta n} (V\phi)(t) \frac{dt}{t}. \notag
\end{align}
Notice that $(V\phi)(t) = \phi([t+1])$. Then, we have
\[\int_n^{\theta n} (V\phi)(t) \frac{dt}{t} = \sum_{i=n}^{[\theta n]-1} \phi(i+1) \int_i^{i+1} \frac{dt}{t} + \int_{[\theta n]}^{\theta n} \phi([\theta n ]+1) \frac{dt}{t}. \]
It is obvious that the second term tends to $0$ as $n \to \infty$. Further, observe that
\begin{align}
&\Bigg|\sum_{i=n}^{[\theta n]-1} \phi(i+1) \int_i^{i+1} \frac{dt}{t} - \sum_{i \in [n, \theta n]} \frac{\phi(i+1)}{i} \Bigg| \notag \\
&\le \Bigg|\phi(n+1) \cdot \left(\frac{1}{n} - \log\left(1+\frac{1}{n}\right)\right) + \phi(n+2) \cdot \left(\frac{1}{n+1} - \log\left(1+\frac{1}{n+1}\right)\right) + \cdots \notag \\
&+ \phi([\theta n]) \cdot \left(\frac{1}{[\theta n]-1} - \log\left(1+\frac{1}{[\theta n]-1}\right)\right)\Bigg| + \left|\frac{\phi([\theta n] + 1)}{[\theta n]}\right| \notag \\
&\le |\phi(n+1)| \cdot \frac{1}{2} \cdot \frac{1}{n^2} + |\phi(n+2)| \cdot \frac{1}{2} \cdot \frac{1}{(n+1)^2} + \cdots \notag \\
&+ |\phi([\theta n])| \cdot \frac{1}{2} \cdot \frac{1}{([\theta n]-1)^2} + \left|\frac{\phi([\theta n] + 1)}{[\theta n]}\right| \notag \\
&\le \frac{\|\phi\|_{\infty}}{2} \left(\frac{1}{n^2} + \frac{1}{(n+1)^2} + \ldots \frac{1}{([\theta n]-1)^2}\right) + \frac{\|\phi\|_{\infty}}{[\theta n]}. \notag
\end{align}
The last expression tends to $0$ as $n \to \infty$. Furthermore, observe that the difference
\[\left|\sum_{i \in [n, \theta n]} \frac{\phi(i)}{i} - \sum_{i \in [n, \theta n]} \frac{\phi(i+1)}{i}\right| \]
tends to $0$ as $n \to \infty$.
Hence, we obtain the following:
\[\overline{p}_C(\phi) = \sup_{\psi \in \mathscr{M}_C} \psi(\phi) = \lim_{\theta \to \infty} \limsup_{n \to \infty} \frac{1}{\log \theta} \sum_{i \in [n, \theta n]} \frac{\phi(i)}{i}. \]

We now consider the functional $\overline{C}_{\infty}$. Note that for $\phi \in l_{\infty}$,
\[C\phi(n) = \frac{1}{n} \sum_{i=1}^n \phi(i) = \frac{1}{n} \int_0^{n} (V\phi)(t)dt = (HV\phi)(n). \]
Then, we have
\[\overline{C}(\phi) = \limsup_{n \to \infty} (C\phi)(n) = \limsup_{n \to \infty} (HV\phi)(n) = \overline{H}(V\phi). \]
Hence, for each $k \ge 2$, we have
\[\overline{C}_k(\phi) = \overline{C}_{k-1}(C\phi) = \overline{H}_{k-1}(HV\phi) = \overline{H}_k(V\phi). \]
Therefore, we obtain
\[\overline{C}_{\infty}(\phi) = \lim_{k \to \infty} \overline{C}_k(\phi) = \lim_{k \to \infty} \overline{H}_k(V\phi) = \overline{Q}(V\phi) = \overline{p}_C(\phi) \]
from the above computation. We have obtained the analytic expression of $\overline{C}_{\infty}(\phi)$:
\[\overline{C}_{\infty}(\phi) =  \lim_{\theta \to \infty} \limsup_{n \to \infty} \frac{1}{\log \theta} \sum_{i \in [n, \theta n]} \frac{\phi(i)}{i}, \quad \phi \in l_{\infty}. \]

The proof of the following characterization of $C_{\infty}$ summability is similar to that of Theorem 6.7. Therefore, the proof is omitted.
\begin{thm}
For $\phi \in l_{\infty}$, $\phi$ is $C_{\infty}$ summable to the number $\alpha$ if and only if
\[\lim_{\theta \to \infty} \frac{1}{\log \theta} \sum_{i \in [n, \theta n]} \frac{\phi(i)}{i} = \alpha \]
uniformly in $n \in \mathbb{N}$.
\end{thm}
It is apparent that this summability method has relation to the logarithmic method, which is defined as follows. For a function $\phi$ on $\mathbb{N}$, we say that $\phi$ is summable to $\alpha$ by the logarithmic method if 
\[\lim_{n \to \infty} \frac{1}{\log n} \sum_{i =1}^n \frac{\phi(i)}{i} = \alpha. \]
The following result, which was given in \cite{Dia}, is an immediate consequence of Theorem 7.1.
\begin{thm}
For $\phi \in l_{\infty}$, if $\phi$ is $C_{\infty}$ summable, then $\phi$ is logarithmic summable.
\end{thm}

\end{document}